\documentstyle{article}

\def\R{{\rm I\!R}}

\reversemarginpar
\catcode`\@=11\relax
\newwrite\@unused
\def\typeout#1{{\let\protect\string\immediate\write\@unused{#1}}}
\def\@nnil{\@nil}
\def\@empty{}
\def\@psdonoop#1\@@#2#3{}
\def\@psdo#1:=#2\do#3{\edef\@psdotmp{#2}\ifx\@psdotmp\@empty \else
    \expandafter\@psdoloop#2,\@nil,\@nil\@@#1{#3}\fi}
\def\@psdoloop#1,#2,#3\@@#4#5{\def#4{#1}\ifx #4\@nnil \else
       #5\def#4{#2}\ifx #4\@nnil \else#5\@ipsdoloop #3\@@#4{#5}\fi\fi}
\def\@ipsdoloop#1,#2\@@#3#4{\def#3{#1}\ifx #3\@nnil
       \let\@nextwhile=\@psdonoop \else
      #4\relax\let\@nextwhile=\@ipsdoloop\fi\@nextwhile#2\@@#3{#4}}
\def\@tpsdo#1:=#2\do#3{\xdef\@psdotmp{#2}\ifx\@psdotmp\@empty \else
    \@tpsdoloop#2\@nil\@nil\@@#1{#3}\fi}
\def\@tpsdoloop#1#2\@@#3#4{\def#3{#1}\ifx #3\@nnil
       \let\@nextwhile=\@psdonoop \else
      #4\relax\let\@nextwhile=\@tpsdoloop\fi\@nextwhile#2\@@#3{#4}}
\def\psdraft{
    \def\@psdraft{0}
}
\def\psfull{
    \def\@psdraft{100}
} \psfull
\newif\if@prologfile
\newif\if@postlogfile
\newif\if@bbllx
\newif\if@bblly
\newif\if@bburx
\newif\if@bbury
\newif\if@height
\newif\if@width
\newif\if@rheight
\newif\if@rwidth
\newif\if@clip
\def\@p@@sclip#1{\@cliptrue}
\def\@p@@sfile#1{
           \def\@p@sfile{#1}
}
\def\@p@@sfigure#1{\def\@p@sfile{#1}}
\def\@p@@sbbllx#1{
        \@bbllxtrue
        \dimen100=#1
        \edef\@p@sbbllx{\number\dimen100}
}
\def\@p@@sbblly#1{
        \@bbllytrue
        \dimen100=#1
        \edef\@p@sbblly{\number\dimen100}
}
\def\@p@@sbburx#1{
        \@bburxtrue
        \dimen100=#1
        \edef\@p@sbburx{\number\dimen100}
}
\def\@p@@sbbury#1{
        \@bburytrue
        \dimen100=#1
        \edef\@p@sbbury{\number\dimen100}
}
\def\@p@@sheight#1{
        \@heighttrue
        \dimen100=#1
        \edef\@p@sheight{\number\dimen100}
}
\def\@p@@swidth#1{
        \@widthtrue
        \dimen100=#1
        \edef\@p@swidth{\number\dimen100}
}
\def\@p@@srheight#1{
        \@rheighttrue
        \dimen100=#1
        \edef\@p@srheight{\number\dimen100}
}
\def\@p@@srwidth#1{
        \@rwidthtrue
        \dimen100=#1
        \edef\@p@srwidth{\number\dimen100}
}
\def\@p@@sprolog#1{\@prologfiletrue\def\@prologfileval{#1}}
\def\@p@@spostlog#1{\@postlogfiletrue\def\@postlogfileval{#1}}
\def\@cs@name#1{\csname #1\endcsname}
\def\@setparms#1=#2,{\@cs@name{@p@@s#1}{#2}}
\def\ps@init@parms{
        \@bbllxfalse \@bbllyfalse
        \@bburxfalse \@bburyfalse
        \@heightfalse \@widthfalse
        \@rheightfalse \@rwidthfalse
        \def\@p@sbbllx{}\def\@p@sbblly{}
        \def\@p@sbburx{}\def\@p@sbbury{}
        \def\@p@sheight{}\def\@p@swidth{}
        \def\@p@srheight{}\def\@p@srwidth{}
        \def\@p@sfile{}
        \def\@p@scost{10}
        \def\@sc{}
        \@prologfilefalse
        \@postlogfilefalse
        \@clipfalse
}
\def\parse@ps@parms#1{
        \@psdo\@psfiga:=#1\do
           {\expandafter\@setparms\@psfiga,}}
\newif\ifno@bb
\newif\ifnot@eof
\newread\ps@stream
\def\bb@missing{
    \typeout{psfig: searching \@p@sfile \space  for bounding box}
    \openin\ps@stream=\@p@sfile
    \no@bbtrue
    \not@eoftrue
    \catcode`\%=12
    \loop
        \read\ps@stream to \line@in
        \global\toks200=\expandafter{\line@in}
        \ifeof\ps@stream \not@eoffalse \fi
        \@bbtest{\toks200}
        \if@bbmatch\not@eoffalse\expandafter\bb@cull\the\toks200\fi
    \ifnot@eof \repeat
    \catcode`\%=14
} \catcode`\%=12
\newif\if@bbmatch
\def\@bbtest#1{\expandafter\@a@\the#1
\long\def\@a@#1
\long\def\bb@cull#1 #2 #3 #4 #5 {
    \dimen100=#2 bp\edef\@p@sbbllx{\number\dimen100}
    \dimen100=#3 bp\edef\@p@sbblly{\number\dimen100}
    \dimen100=#4 bp\edef\@p@sbburx{\number\dimen100}
    \dimen100=#5 bp\edef\@p@sbbury{\number\dimen100}
    \no@bbfalse
} \catcode`\%=14
\def\compute@bb{
        \no@bbfalse
        \if@bbllx \else \no@bbtrue \fi
        \if@bblly \else \no@bbtrue \fi
        \if@bburx \else \no@bbtrue \fi
        \if@bbury \else \no@bbtrue \fi
        \ifno@bb \bb@missing \fi
        \ifno@bb \typeout{FATAL ERROR: no bb supplied or found}
            \no-bb-error
        \fi
        \count203=\@p@sbburx
        \count204=\@p@sbbury
        \advance\count203 by -\@p@sbbllx
        \advance\count204 by -\@p@sbblly
        \edef\@bbw{\number\count203}
        \edef\@bbh{\number\count204}
}
\def\in@hundreds#1#2#3{\count240=#2 \count241=#3
             \count100=\count240    
             \divide\count100 by \count241
             \count101=\count100
             \multiply\count101 by \count241
             \advance\count240 by -\count101
             \multiply\count240 by 10
             \count101=\count240    
             \divide\count101 by \count241
             \count102=\count101
             \multiply\count102 by \count241
             \advance\count240 by -\count102
             \multiply\count240 by 10
             \count102=\count240    
             \divide\count102 by \count241
             \count200=#1\count205=0
             \count201=\count200
            \multiply\count201 by \count100
            \advance\count205 by \count201
             \count201=\count200
            \divide\count201 by 10
            \multiply\count201 by \count101
            \advance\count205 by \count201
             \count201=\count200
            \divide\count201 by 100
            \multiply\count201 by \count102
            \advance\count205 by \count201
             \edef\@result{\number\count205}
}
\def\compute@wfromh{
        \in@hundreds{\@p@sheight}{\@bbw}{\@bbh}
        \edef\@p@swidth{\@result}
}
\def\compute@hfromw{
        \in@hundreds{\@p@swidth}{\@bbh}{\@bbw}
        \edef\@p@sheight{\@result}
}
\def\compute@handw{
        \if@height
            \if@width
            \else
                \compute@wfromh
            \fi
        \else
            \if@width
                \compute@hfromw
            \else
                \edef\@p@sheight{\@bbh}
                \edef\@p@swidth{\@bbw}
            \fi
        \fi
}
\def\compute@resv{
        \if@rheight \else \edef\@p@srheight{\@p@sheight} \fi
        \if@rwidth \else \edef\@p@srwidth{\@p@swidth} \fi
}
\def\compute@sizes{
    \compute@bb
    \compute@handw
    \compute@resv
}
\def\psfig#1{\vbox {
    \ps@init@parms
    \parse@ps@parms{#1}
    \compute@sizes
    \ifnum\@p@scost<\@psdraft{
        \typeout{psfig: including \@p@sfile \space }
        \special{ps::[begin]    \@p@swidth \space \@p@sheight \space
                \@p@sbbllx \space \@p@sbblly \space
                \@p@sbburx \space \@p@sbbury \space
                startTexFig \space }
        \if@clip{
            \typeout{(clip)}
            \special{ps:: \@p@sbbllx \space \@p@sbblly \space
                \@p@sbburx \space \@p@sbbury \space
                doclip \space }
        }\fi
        \if@prologfile
            \special{ps: plotfile \@prologfileval \space } \fi
        \special{ps: plotfile \@p@sfile \space }
        \if@postlogfile
            \special{ps: plotfile \@postlogfileval \space } \fi
        \special{ps::[end] endTexFig \space }
        \vbox to \@p@srheight true sp{
            \hbox to \@p@srwidth true sp{
                \hfil
            }
        \vfil
        }
    }\else{
        \vbox to \@p@srheight true sp{
        \vss
            \hbox to \@p@srwidth true sp{
                \hss
                \@p@sfile
                \hss
            }
        \vss
        }
    }\fi
}} \catcode`\@=12\relax

\newenvironment{pf}{\noindent {\bf Proof:}}{{\hfill $\Box$}}
\newenvironment{proofT}{\noindent {\bf Proof of Theorem}}{{\hfill $\Box$}}

\def\noproof{\unskip{\hfill $\Box$}}

\def\qed{$\Box$}

\newcommand{\lemlab}[1]{\label{lemma:#1}}
\newcommand{\thmlab}[1]{\label{thm:#1}}
\newcommand{\codelab}[1]{\label{code:#1}}
\newcommand{\alglab}[1]{\label{alg:#1}}
\newcommand{\eqlab}[1]{\label{eq:#1}}
\newcommand{\corlab}[1]{\label{cor:#1}}
\newcommand{\deflab}[1]{\label{def:#1}}
\newcommand{\tablab}[1]{\label{tab:#1}}
\newcommand{\figlab}[1]{\label{fig:#1}}
\newcommand{\seclab}[1]{\label{section:#1}}
\newcommand{\chaplab}[1]{\label{chapter:#1}}

\newcommand{\lemref}[1]{\ref{lemma:#1}}
\newcommand{\coderef}[1]{\ref{code:#1}}
\newcommand{\thmref}[1]{\ref{thm:#1}}
\newcommand{\algref}[1]{\ref{alg:#1}}
\newcommand{\corref}[1]{\ref{cor:#1}}
\newcommand{\defref}[1]{\ref{def:#1}}
\newcommand{\tabref}[1]{\ref{tab:#1}}
\newcommand{\figref}[1]{\ref{fig:#1}}
\newcommand{\eqref}[1]{(\ref{eq:#1})}
\newcommand{\secref}[1]{\ref{section:#1}}
\newcommand{\chapref}[1]{\ref{chapter:#1}}

\newcommand{\li}{\item}
\newtheorem{theorem}{Theorem}[section]
\newtheorem{lemma}[theorem]{Lemma}

\newtheorem{cor}[theorem]{Corollary}
\newtheorem{conj}{Conjecture}

\newtheorem{open}[conj]{Open Problem}

{\catcode`\@=11
\gdef\setft#1#2#3{%
\def\@oddfoot{
{\setbox0=\hbox{#1} \setbox1=\hbox{#3} \ifdim\wd0>\wd1
\dimen0=\wd0 \box0\hfil#2\hfil\hbox to\dimen0{\hfil\hfil\box1}
\else \dimen0=\wd1 \hbox to\dimen0{\box0\hfil }\hfil#2\hfil\box1
\fi }}} }

\def\complaint#1{}
\def\withcomplaints{
\newcounter{mycomplaints}
\def\complaint##1{\refstepcounter{mycomplaints}%
\ifhmode%
\unskip%
{\dimen1=\baselineskip \divide\dimen1 by 2 %
\raise\dimen1\llap{\tiny -\themycomplaints-}}\fi%
\marginpar{\tiny [\themycomplaints]: ##1}}%
}

\withcomplaints 

\begin{document}

\title{Planar Minimally Rigid Graphs and
Pseudo-Triangulations}

\author{
\setcounter{footnote}{0}%
\def\thefootnote{\arabic{footnote}}
Ruth Haas \footnotemark[1] \hspace{3mm}%
David Orden \footnotemark[2] \hspace{3mm}%
G\"unter Rote \footnotemark[3]\hspace{3mm}%
\\%
\ \setcounter{footnote}{3} %
\def\thefootnote{\arabic{footnote}}%
Francisco Santos \footnotemark[2] \hspace{3mm}%
Brigitte Servatius \footnotemark[4] \hspace{3mm}%
Herman Servatius \footnotemark[4]\hspace{3mm}%
\\%
\setcounter{footnote}{6}%
\def\thefootnote{\arabic{footnote}}%
Diane ~Souvaine \footnotemark[5]\hspace{3mm} %
Ileana ~Streinu \footnotemark[6]\hspace{3mm}%
and Walter ~Whiteley \footnotemark[7]\hspace{3mm}%
}

\date{}

{
\footnotetext[1]{Department\ of\ Mathematics, Smith College,
Northampton, MA 01063, USA, rhaas@math.smith.edu.}%
\footnotetext[2]{Departamento de Matematicas, Estadistica y
Computacion, Universidad de Cantabria, E-39005 Santander, Spain,
{\{ordend, santos\}}@matesco.unican.es. Supported by grant
BFM2001-1153 of the Spanish Ministry of Science and Technology.}%
\footnotetext[3]{Institut f\"ur Informatik, Freie Universit\"at
Berlin, Takustra{\ss}e 9, D-14195~Berlin,~Germany.
\textsl{rote@inf.fu-berlin.de}. Partly supported by the Deutsche
Forschungsgemeinschaft (DFG) under grant RO~2338/2-1.}%
\footnotetext[4]{Mathematics Department, Worcester Polytechnic
Institute, Worcester, MA 01609,
USA. \textsl{{\{bservat, hservat\}}@math.wpi.edu}.}%
\footnotetext[5]{Department of Computer Science, Tufts University,
Medford MA, USA. \textsl{dls@cs.tufts.edu}. Supported by NSF Grant
EIA-9996237.}%
\footnotetext[6]{Department of Computer Science, Smith College,
Northampton, MA 01063, USA. \textsl{streinu@cs.smith.edu}.
Supported by NSF grants CCR-0105507 and CCR-0138374.}%
\footnotetext[7]{Department of Mathematics and Statistics, York
University, Toronto, Canada. \textsl{whiteley@mathstat.yorku.ca}.
Supported by NSERC (Canada) and NIH(USA).}%
}

\maketitle

\begin{abstract}
Pointed pseudo-triangulations are planar minimally rigid graphs
embedded in the plane with {\em pointed} vertices (adjacent to an
angle larger than $\pi$). In this paper we prove that the opposite
statement is also true, namely that planar minimally rigid graphs
always admit pointed embeddings, even under certain natural
topological and combinatorial constraints. The proofs yield
efficient embedding algorithms. They also provide - to the best of
our knowledge - the first algorithmically effective result on
graph embeddings with oriented matroid constraints other than
convexity of faces. These constraints are described by {\em
combinatorial pseudo-triangulations}, first defined and studied in
this paper. Also of interest are our two proof techniques, one
based on Henneberg inductive constructions from combinatorial
rigidity theory, the other on a generalization of Tutte's
barycentric embeddings to directed graphs.

\end{abstract}


\pagestyle{plain}

\section{Introduction}

\label{introduction}

In this paper we bring together two classical topics in graph
theory, planarity and rigidity, to answer the fundamental question
(posed in \cite{streinu2}) of {\it characterizing the class of
planar graphs which admit pointed pseudo-triangular embeddings}.
Our main result is that this coincides with the class of all {\it
planar minimally rigid graphs (planar Laman graphs)}. Furthermore
we extend the result in several directions, attacking the same
type of question for other (not necessarily pointed) classes of
pseudo-triangulations and for {\it combinatorial
pseudo-triangulations}, a new class of objects first introduced
and studied in this paper.

\medskip

\noindent {\bf Novelty.} As opposed to traditional planar graph
embeddings, where all the faces are designed to be convex, ours
have interior faces which are as non-convex as possible
(pseudo-triangles). Planar graph embeddings with {\it non-convex}
faces have not been systematically studied. Our result links them
to rigidity theoretic and matroidal properties of planar graphs.
To the best of our knowledge, this is the first result holding for
an interesting family of graphs on algorithmically efficient graph
embeddings with oriented matroid constraints other than convexity
of faces. In contrast, the universality theorem for pseudo-line
arrangements of Mn\"ev \cite{mnev} implies that the general
problem of embedding graphs with oriented matroid constraints is
as hard as the existential theory of the reals.

\medskip
\noindent {\bf Proof Techniques and Algorithmic Results.} We
present two proof techniques of independent interest. The first
one is of a {\it local} nature, relying on incremental (inductive)
constructions known in rigidity theory as {\it Henneberg
constructions}. The second one is based on a {\it global} approach
making use of a directed version of {\it Tutte's barycentric
embeddings}, which - to the best of our knowledge - is first
proven in this paper. Both proofs are constructive, yield
efficient algorithms, emphasize distinct aspects of the result and
lead into new directions of further investigation: combinatorial
versus geometric embeddings, local versus global coordinate
finding.

\medskip
\noindent {\bf Laman Graphs and Pseudo-Triangulations.} Let
$G=(V,E)$ be a graph with $n=|V|$ vertices and $m=|E|$ edges. $G$
is a {\it Laman graph} if $m=2n-3$ and every subset of $k$
vertices spans at most $2k-3$ edges.\footnote{This is called the
definition {\it by counts} of Laman graphs. An equivalent
definition via Henneberg constructions will be given in
Section~\ref{section:preliminaries}.} An embedding $G(P)$ of the
graph $G$ on a set of points $P=\{ p_1, \cdots, p_n\}\subset R^2$
is a mapping of the vertices $V$ to points in the Euclidean plane
$i\mapsto p_i\in P$. The edges $ij\in E$ are mapped to straight
line segments $p_ip_j$. We say that the vertex $i$ of the
embedding $G(P)$ is {\it pointed} if all its adjacent edges lie on
one side of some line through $p_i$. Equivalently, some
consecutive pair of edges adjacent to $i$ (in the circular
counter-clockwise order around the vertex) spans a {\it reflex}
angle. An embedding $G(P)$ is {\it non-crossing} if no pair of
segments $p_i p_j$ and $p_kp_l$ corresponding to non-adjacent
edges $ij, kl\in E, i,j\not\in\{k,l\}$ have a point in common. A
graph $G$ is {\it planar} if it has a non-crossing embedding.

A {\it pseudo-triangle} is a simple planar polygon with exactly
three convex vertices. A {\it pseudo-triangulation} of a planar
set of points is a non-crossing embedded graph $G(P)$ whose outer
face is convex and all interior faces are pseudo-triangles. In a
{\it pointed pseudo-triangulation} all the vertices are pointed.
See Figure~\figref{pseudoT}.

\begin{figure}[ht]
\begin{center}
\ \psfig{figure=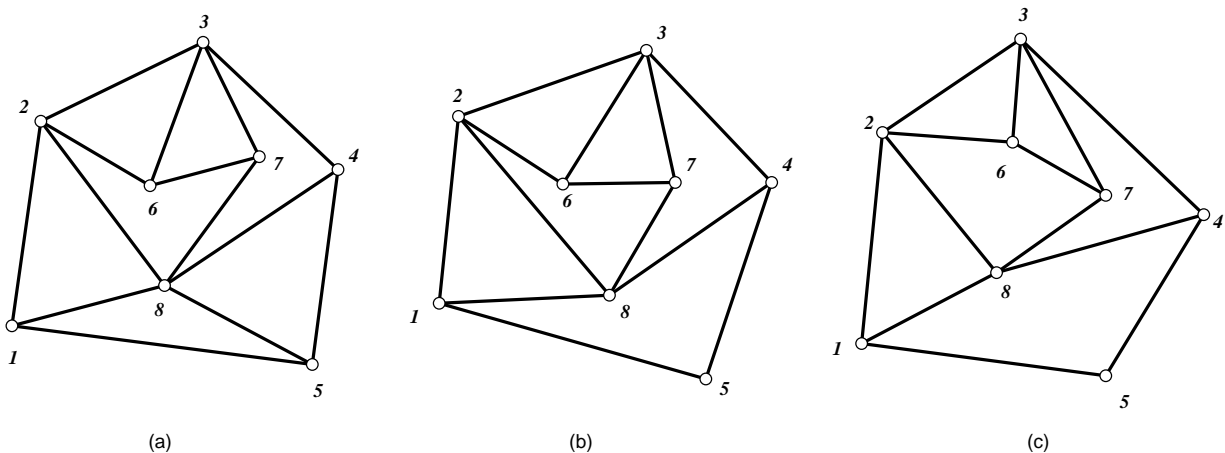,height=2.7in}
\end{center}
\caption{(a) A pseudo-triangulation (necessarily non-pointed,
since the underlying graph is a circuit, not a Laman graph) and
two embeddings of a planar Laman graph: (b) is a pointed
pseudo-triangulation, (c) is not: the faces $2876$ and $1548$ are
not pseudo-triangles and the vertices $6$ and $8$ are not
pointed.} \figlab{pseudoT}
\end{figure}

Pseudo-triangulations are thus planar graphs with a special
embedding. As shown in \cite{streinu}, the graphs of pointed
pseudo-triangulations are Laman graphs. They have very useful
rigidity theoretic properties and a wealth of combinatorial
properties. They have emerged via several applications in
Computational Geometry, where they are used in Kinetic Data
Structures for collision detection and in certain motion planning
problems.

\medskip
\noindent {\bf Historical Perspective.} Techniques from Rigidity
Theory have been recently applied to problems such as collision
free robot arm motion planning \cite{cdr, streinu}, molecular
conformations \cite{JRKT, whiteleyM,thorpe} or sensor and network
topologies with distance and angle constraints \cite{eren}.

Laman graphs are {\it the} fundamental objects in $2$-dimensional
Rigidity Theory. Also known as {\it isostatic} or {\it generically
minimally rigid graphs}, they characterize combinatorially the
property that a graph, embedded on a generic set of points in the
plane, is infinitesimally rigid (with respect to the induced edge
lengths). See \cite{laman}, \cite{gss}, \cite{walter}. The most
famous open question in Rigidity Theory (the {\it Rigidity
Conjecture}, see \cite{gss}) is finding their $3$-dimensional
counterpart.

Pseudo-triangulations are relatively new objects, introduced and
applied in Computational Geometry for problems such as visibility
\cite{pocchiola}, \cite{pocchiola1}, \cite{st}, kinetic data
structures \cite{guibas} and motion planning for robot arms
\cite{streinu}. They have rich combinatorial, rigidity theoretic
and polyhedral properties (\cite{streinu}, \cite{rss},
\cite{os-pncgp-02}), many of which have only recently started to
be investigated (\cite{rrss}, \cite{kettner}, \cite{ass},
\cite{bespat}, \cite{aaks}, \cite{bron}, \cite{AOSS03}). In
particular, the fact that they are Laman graphs which become {\it
expansive} mechanisms when one of their convex hull edges is
removed, has proven to be crucial in designing efficient motion
planning algorithms for planar robot arms, see \cite{streinu}.
Finding their $3$-dimensional counterpart, which is perhaps the
main open question about pseudo-triangulations and expansive
motions, may lead to efficient motion planning algorithms for
certain classes of $3$-dimensional linkages, with potential impact
on understanding protein folding processes.

Graph Drawing is a field with a distinguished history, and
embeddings of planar graphs have received substantial attention in
the literature (\cite{fary},\cite{tutteC}, \cite{tutteG},
\cite{cvpv}, \cite{fpp}, \cite{schnyder}, \cite{dett}). Extensions
of graph embeddings from straight-line to pseudo-line
segments have been recently considered (see e.g. \cite{pt}). It is
natural to ask which such embeddings are stretchable, i.e. whether
they can be realized with straight-line segments while maintaining
some desired combinatorial substructure. Indeed, the primordial
planar graph embedding result, Fary's Theorem \cite{fary}, is just
an instance of answering such a question. Graph embedding
stretchability questions have usually ignored oriented matroidal
constraints, allowing for the free reorientation of triplets of
points when not violating other combinatorial conditions. The
notable exception concerns the still widely open {\it visibility
graph recognition problem}, approached in the context of
pseudo-line arrangements (oriented matroids) by \cite{ors}. In
\cite{starlike} it is shown that it is not always possible to
realize with straight-lines a pseudo-visibility graph, while
maintaining oriented matroidal constraints.

In contrast,  {\it this paper gives the first non-trivial
stretchability result on a natural graph embedding problem with
oriented matroid constraints.} It adds to the already rich body of
surprisingly simple and elegant combinatorial properties of
pointed pseudo-triangulations by proving a natural connection.

\medskip
\noindent {\bf Main Result.} We are interested in {\it planar
Laman graphs}. Not all Laman graphs fall into this category. For
example, $K_{3,3}$ is Laman but not planar. But the underlying
graphs of all pointed pseudo-triangulations are planar Laman.
We prove that the converse is always true:

\begin{theorem}

{\bf (Main Theorem)}
\label{mainThm}

Every planar Laman graph can be embedded as a pointed
pseudo-triangulation.
\end{theorem}

\noindent The following characterization follows then from well
known properties of Laman graphs:

\begin{cor} \label{cor:equivalences-lamangraphs}
Given a {\em planar} graph $G$, the following conditions are
equivalent:
\begin{itemize}
\item[(i)] $G$ is a Laman graph%
\item[(ii)] Generically, $G$ is minimally rigid%
\item[(iii)] $G$ can be embedded as a pointed pseudo-triangulation%
\end{itemize}
\end{cor}

We prove in fact several stronger results, allowing the {\em a
priori} choice of the facial structure (Theorem \ref{thm:mainL})
and even of the combinatorial information regarding which vertices
are convex in each face (Theorem \ref{thm:mainC}). This last
result needs the apparatus of {\em combinatorial
pseudo-triangulations}, first defined and studied in this paper.

Finally, we answer a natural question related to the underlying
matroidal structure of planar rigidity and extend the result to
planar rigidity circuits, which are minimal dependent sets in the
rigidity matroid where the bases (maximally independent sets) are
the Laman graphs. By adding edges to a pointed (minimum)
pseudo-triangulation while maintaining planarity, the graph has
increased dependency level (in the rigidity matroid) and can no
longer be realized with all vertices pointed, but it can always be
realized with straight edges. Our concern is to maintain the
minimum number of non-pointed vertices, for the given edge count.
For circuits, this number is one, and we show that it can be
attained.

\medskip
\noindent {\bf Organization.} The paper is organized as follows.
In Section~\ref{section:preliminaries} we give the basic
definitions needed for an independent reading of
Section~\ref{section:main}. For increased readability, additional
technical definitions are later included in the sections that use
them. The first proof of the Main Theorem is presented in
Section~\ref{section:main}, which is further devoted to all the
proofs (combinatorial or geometric) making use of the inductive
Henneberg technique: planar Laman graphs, combinatorial
pseudo-triangulations, pointed pseudo-triangulations and
Laman-plus-one combinatorial and geometric pseudo-triangulations.
Section~\ref{section:cpt} is devoted to combinatorial
pseudo-triangulations and to the perfect matching technique for
assigning combinatorial pseudo-triangular labelings to plane
graphs. Section~\ref{section:proofT} focuses on the second proof
technique based on Tutte embeddings and contains our most general
result on plane graph embeddings compatible with {\it given}
combinatorial pseudo-triangulations. We conclude in
Section~\ref{section:open} with a list of further directions of
research and open questions.

\section{Preliminaries}
\label{section:preliminaries}

For the standard graph and rigidity theoretical terminology used
in this paper we refer the reader to \cite{gss} and
\cite{walterS}. For relevant facts about pointed
pseudo-triangulations, see \cite{streinu}. In this section we
continue what we started in the Introduction and give most of the
definitions needed for reading Sections \ref{section:main} and
\ref{section:cpt}. The technically denser Section
\ref{section:proofT} contains its own additional concepts.

\medskip
\noindent {\bf Notation and Abbreviations.} Throughout the paper
we will abbreviate {\it counter-clockwise} as {\it ccw}.  To
emphasize that a graph has $n$ vertices we may denote it by $G_n$.
We will occasionally abbreviate {\it combinatorial
pseudo-triangulation} as {\it cpt} and {\it pointed combinatorial
pseudo-triangulation} as {\it pointed cpt}.

\medskip
\noindent {\bf Plane Graphs.} A non-crossing embedding of a
connected planar graph $G$ partitions the plane into faces
(bounded or unbounded), edges 
and vertices.
Their incidences are fully captured by the vertex {\it rotations}:
the ccw circular order of the edges incident to each vertex in the
embedding. A {\it sphere embedding} of a planar graph refers to a
choice of a system of rotations (and thus of a facial structure),
and is oblivious of an outer face.
It is well-known (Whitney \cite{Whitney}) that a $3$-connected
planar graph induces a unique set of rotations,
but $2$-connected ones may induce several. A {\it plane graph} is
a spherical graph with a choice of a particular face as the outer
face. Every simple plane graph can be realized with straight-line
edges in the plane (Fary's theorem \cite{fary}).

A (combinatorial) {\it angle} (incident to a vertex or a face in a
plane graph) is a pair of consecutive  edges (consecutive in the
order given by the rotations) incident to the vertex or face.


\medskip
\noindent {\bf Pseudo-Triangulations.} We have defined
pseudo-triangles, pseudo-triangulations and pointed
pseudo-triangulations in the Introduction. In addition, we will
use the following related concepts. The {\it corners} of a
pseudo-triangle are its three convex angles, and its {\it side
chains} are the pieces of the boundary between two corners
(vertices and edges). The {\it extreme edges} of a pointed vertex
are the two edges incident with its unique incident reflex angle.

\medskip
\noindent {\bf Minimally Rigid (Laman) Graphs and Henneberg
constructions.} Besides the definition by {\it counts} given in
the Introduction, Laman graphs can be characterized in a variety
of ways. In particular, a Laman graph on $n$ vertices has an
inductive construction as follows (see \cite{henneberg, walter}).
Start with an edge for $n=2$. At each step, add a new vertex in
one of the following two ways:

\begin{itemize}
\item {\bf Henneberg I} (vertex addition): the new vertex is connected
via two new edges to two old vertices.
\item {\bf Henneberg II} (edge splitting): a new vertex
is added on some edge (thus splitting the edge into two new edges)
and then connected to a third vertex. Equivalently, this can be
seen as removing an edge, then adding a new vertex connected to
its two endpoints and to some other vertex.
\end{itemize}

See Figure~\figref{henne}, where we show drawings with crossing
edges, to emphasize that the Henneberg constructions work for {\it
general}, not necessarily {\it planar} Laman graphs.

\begin{figure}[ht]
\begin{center}
\ \psfig{figure=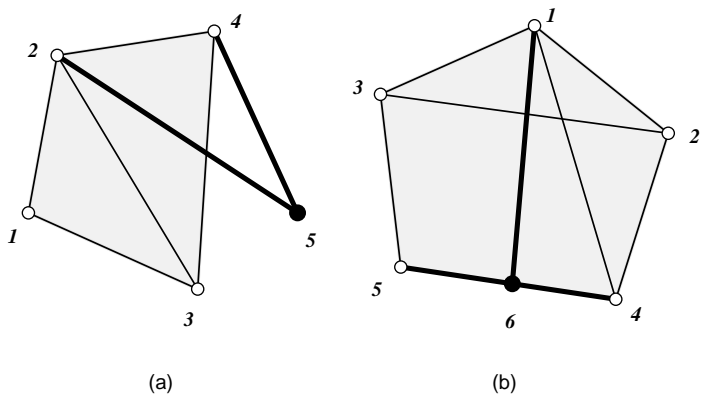,height=2.2in}
\end{center}
\caption{
Illustration of the two types of steps in a Henneberg sequence,
with vertices labelled in the construction order. The shaded part
is the old graph, to which the black vertex is added. (a)
Henneberg I for vertex $5$, connected to old vertices $3$ and $4$.
(b) Henneberg II for vertex $6$, connected to old vertices $3$,
$4$ and $5$. }%
\figlab{henne}
\end{figure}

We will make heavy use of the following result, essentially stated
by Henneberg \cite{henneberg}, and of its proof, due to Tay and
Whiteley \cite{taywhiteley}.

\begin{lemma}
\label{lemma:henneberg} A graph is Laman if and only if it has a
Henneberg construction.
\end{lemma}


\smallskip

The proof of Lemma \ref{lemma:henneberg} proceeds inductively to
show that there always exists a vertex of degree $2$ or $3$ which
can be removed in the reverse order of a Henneberg step while
maintaining the Laman property. It is instructive to give a
slightly more general proof. We will make use of it in
Section~\ref{section:main}.

\begin{lemma}
\label{lemma:hennebergS} A Laman graph has a Henneberg
construction starting from any prescribed subset of two vertices.
Moreover, if there exist three vertices of degree $3$ mutually
connected in a triangle, then we can prescribe them as the first
three vertices of the Henneberg construction.
\end{lemma}

\begin{pf} Let $G_n=(V,E)$ be a Laman graph on $n=|V|$ vertices and
let $V_2\subset V$ be any subset of two vertices. We show that as
long as $n > 2$ we can always remove a vertex not in $V_2$ in the
opposite direction of a Henneberg step. In the actual Henneberg
construction this amounts to starting the induction from this
prescribed pair.

Vertices of degree $0$ or $1$ do not exist in Laman graphs,
otherwise the Laman property would be violated on a subset of
$n-1$ vertices. Since $G_n$ has $2n-3$ edges, a simple count shows
that there exist at least three vertices of degree at most $3$,
hence either of degree $2$ or $3$. At least one of them (call it
$v$) is not in $V_2$. This will be the vertex we choose to remove,
in a backwards application of a Henneberg step.

If $v$ has degree $2$, we remove $v$ and its incident edges: the
resulting graph on $n-1$ vertices and $2(n-1)-1$ edges is clearly
Laman. If $v$ has degree $3$, let its neighbors be $v_1$, $v_2$
and $v_3$. The removal of $v$ and of its three adjacent edges
produces a graph $G'_{n-1}$ with a deficit of one edge: $n-1$
vertices but only $2(n-1)-4$ edges. We must put back one edge
joining one of the three pairs of vertices in $v_1, v_2, v_3$.
Consider the {\em rigid components} of $G'_{n-1}$: maximal subsets
of some $k$ vertices spanning $2k-3$ edges.
The three endpoints $v_1, v_2$ and $v_3$ cannot belong to the same
rigid component (otherwise the Laman count would be violated  in
$G_n$  on the subset consisting of this component and $v$). Two
rigid components share at most one vertex, otherwise their union
would be a larger Laman subgraph. Suppose that $v_1$ and $v_2$ are
in distinct rigid components. Then adding an edge between $v_1$
and $v_2$ doesn't violate the Laman condition on any subset and
completes $G'_{n-1}$ to a Laman graph $G_{n-1}$.

If $G_n$ contains a subset $V_3$ of three degree $3$ vertices
connected in a triangle, then a similar counting argument shows
that there is an additional vertex of degree at most $3$. This
fourth vertex can be removed in such a way that the invariant (of
having three vertices of degree $3$ connected in a triangle) is
maintained. Hence the vertices of $V_3$ may be prescribed as the
three starting vertices.
To finish, let's show that the invariant is maintained. Let $G'$
be the subgraph induced on the vertices in $V\setminus V_3$: if it
contains more than two elements, then it spans $2n-3-6=2(n-3)-3$
edges, hence it is Laman. Let $N(V_3)$ be the neighbors of $V_3$
in $G_n$: none of these vertices can be of degree $2$, otherwise
$G'$ would not be Laman. If there is a vertex $v'$ of degree $3$
in $N(V_3)$ which is removed at some Henneberg step, one must put
back an edge incident to two of its neighbors, and one of them
must be in $V_3$: otherwise, the induced subgraph $G'$ on
$V\setminus V_3$ (after performing the reverse Henneberg step)
would violate the Laman counts. Since at any reverse Henneberg
step we remove either vertices of degree $2$ (which are not
incident to $V_3$) or of degree $3$, which do not change the
degree of their neighbors in $V_3$, it follows that the vertices
in $V_3$ maintain their degrees and the property of being
connected in a triangle throughout the construction (in fact,
until $n=5$, from which point it is easy to see that $V_3$ can
still be prescribed).
\end{pf}

\medskip
\noindent {\bf Laman-plus-one Graphs and Rigidity Circuits.} A
{\it Laman-plus-one} graph is a Laman graph with one additional
edge. It has $2n-2$ edges and every subset of $k$ edges induces at
most $2k-2$ edges. A {\it rigidity circuit} (shortly, a {\it
circuit}) is a graph with the property that removing {\it any}
edge produces a Laman graph. It is therefore a special
Laman-plus-one graph. In a rigidity circuit $G$ with $n$ vertices,
the number $n$ of vertices is at least $4$, the number $m$ of
edges is $2n-2$ and every subset of $k$ vertices spans at most
$2k-3$ edges. Moreover, the minimum degree in a circuit is $3$. It
is straightforward to prove that a Laman-plus-one graph contains a
unique rigidity circuit: take the maximal subgraph satisfying the
circuit counts. It is unique, because otherwise the union of two
circuits would violate the Laman-plus-one counts.

These concepts are motivated by the matroid view of Rigidity
Theory, see \cite{gss}: Laman graphs correspond to bases (maximal
{\it independent} sets of edges) in the generic rigidity matroid,
while the circuits are the minimally {\it dependent} sets.

It has been proven recently that $3$-connected rigidity circuits
admit an inductive (Henneberg-type) construction (using only
Henneberg II steps and starting from $K_4$), where all
intermediate graphs are themselves circuits. All rigidity circuits
are $2$-connected, hence they can be obtained by making use of
Tutte's Theorem on the structure of $3$-connected graphs in terms
of $2$-connected components (see \cite{Tutte1}, \cite{Tutte}).

We show now (and use later) that Laman-plus-one graphs also admit
a simple Henneberg construction. This proof structure (much easier
than \cite{bergjordan} because of a simpler inductive invariant)
will be used in Section~\ref{section:main} to show stretchability
of planar Laman-plus-one graphs, and thus of planar rigidity
circuits.

\begin{lemma}
\label{lemma:hennebergC} A graph $G$ is Laman-plus-one if and only
if it has a Henneberg construction starting from a $K_4$.
\end{lemma}

\begin{pf}
The proof is similar to that of Lemma \ref{lemma:henneberg} and
uses the $2n-2$ counts. There must be at least two vertices in $G$
of degree at most $3$. Vertices of degree $2$ are outside the
circuit and are removed just as in the Laman case. We now show how
to handle a vertex $v$ of degree $3$. If there are no vertices of
degree $2$, there are at least four of degree $3$.

Let $C=(V_c,E_c)$ be {\em the} unique induced subgraph which is
the circuit of $G$, $V_o=V\setminus V_c$ the vertices outside the
circuit and $V_b\subset V_c$ the {\it boundary} of $V_c$, i.e. the
set of vertices in  $V_c$ incident to a vertex in $V_o$. It is
easy to see (from the $2n-2$ counts) that between $V_b$ and $V_o$
there must be at least two edges if $|V_0|=1$ or three edges if
$|V_o|\geq 2$.

If $v\in V_o$ then its neighbors cannot all three belong to the
circuit, because otherwise the subgraph induced on $V_c\cup \{v\}$
would violate the $2n-2$ counts. Remove (temporarily) an edge $ab$
of $G$ from inside $C$: the resulting graph is Laman, containing
$C$ without this edge as a rigid block (subset on which the Laman
count is satisfied with equality). By Lemma \ref{lemma:hennebergS}
there is a well defined way of removing $v$ and placing back an
edge to perform a Henneberg II step in reverse: the added edge is
not between two vertices of $V_c$. Therefore we can put back the
temporarily removed edge $ab$ to get a Laman-plus-one graph.

If $v\in V_c$, notice first that it cannot be on the boundary
$V_b$, otherwise its degree in $C$ would be at most $2$,
contradicting the fact that $C$ is a circuit. So all its three
neighbors $v_1, v_2$ and $v_3$ are in $V_c$. Removing $v$ and its
incident edges produces a Laman graph. Either all of the edges
$v_1v_2, v_1v_3, v_2v_3$ are present in $G$ or not. If not (say,
$v_1v_2$ is missing), then we add $v_1v_2$, get a Laman-plus-one
graph and continue the induction. Otherwise, the circuit was a
$K_4$. If there are no vertices outside the circuit, we are done.
Otherwise, there are at least two edges out of $V_b$, increasing
the degree of at least one vertex in the circuit $K_4$. If there
are no vertices of degree $2$, then there must be at least one
other vertex of degree $3$ in $V_o$. We will perform the Henneberg
step on it (and thus not touch $K_4$ until the end).
\end{pf}

\medskip
\noindent {\bf Combinatorial Pseudo-Triangulations.}  Let $G$ be a
plane connected graph \footnote{The definition is valid in a more
general setting than what we use in this paper, and works even
with multiple edges, vertices of degree one and loops. In such a
case, a vertex of degree one is  incident to a unique  angle,
labeled {\em big}.}. A {\it combinatorial pseudo-triangulation
(cpt)} of $G$ is an assignment of labels {\it big} (or {\it
reflex}) and {\it small} (or {\it convex}) to the angles of $G$
such that:
\begin{itemize}
\item[(i)] Every face except the outer face gets {\em three}
vertices marked {\it small}. These will be called the {\em corners} of the face.%
\item[(ii)] The outer face gets only {\it big} labels (has no corners).%
\item[(iii)] Each vertex is incident to at most one angle labeled
{\it big}. If it is incident to a big angle, it is called {\em pointed}. %
\item[(iv)] A vertex of degree $2$ is incident to one angle
labeled {\it big}.
\end{itemize}

By analogy with pseudo-triangulations, we also define extreme
edges, side-chains and non-pointed vertices of combinatorial
pseudo-triangulations.

Combinatorial pseudo-triangulations share many combinatorial
properties with pseudo-triangulations. The following lemma follows
easily from the definition.

\begin{lemma}
\label{lem:cptcount}%
A combinatorial pseudo-triangulation on $n$ vertices has at least
$2n-3$ edges. If a cpt has $m\geq 2n-3$ edges, then it contains
exactly $m-(2n-3)$ non-pointed vertices.
\end{lemma}

\begin{pf} Let
$V_K$ be the set of non-pointed vertices. Let $m$ be the number of
edges, $f$ the number of faces, $k$ the size of $V_K$ and $d_v$
the degree of a vertex $v$. We count the number of small angles in
two ways and apply Euler's formula. Summing over the faces we get
$3(f-1)$. Summing over vertices we get $\sum_{v\not\in V_K}(d_v -
1) + \sum_{v\in V_K}d_v= \sum_{v\in V}d_v - (n-k) = 2m - n + k$.
This solves to $m = 2n - 3 + k$ and proves the statement.
\end{pf}

\smallskip

A cpt with exactly $2n-3$ edges will have all vertices pointed:
we'll call it a {\it pointed combinatorial pseudo-triangulation}
(pointed cpt). Another case of interest in this paper is when
exactly one vertex is {\it combinatorially non-pointed}, i.e. has
no incident big angle: we call it a {\it pointed-plus-one} cpt.
Notice that in this case the non-pointed vertex has degree at
least $3$, is interior, i.e. not incident to the outer face and
that the cpt has $2n-2$ edges.

In Section~\ref{section:main} we will prove that all planar Laman
graphs and all planar Laman-plus-one graphs have cpt assignments.
The previous lemma implies that such cpt assignments must be
pointed for Laman graphs, resp. pointed-plus-one for
Laman-plus-one graphs (defined below).


\medskip
\noindent {\bf Pointed-plus-one and Circuit
Pseudo-triangulations.} A {\it pointed-plus-one
pseudo-triangulation} is a pseudo-triangulation with exactly one
non-pointed vertex. It is easy to see that it has $2n-2$ edges,
and is in fact just a planar Laman-plus-one graph embedded as a
pseudo-triangulation. A {\it pseudo-triangulation circuit} is a
planar rigidity circuit embedded as a pseudo-triangulation. A
reminder that by a {\it pseudo-triangulation} we mean any
decomposition into pseudo-triangles, which may not be pointed. In
fact, any such pseudo-triangulation with more than $2n-3$ edges is
necessarily {\it non-pointed}\footnote{In general, those with $k$
non-pointed vertices (all interior) have $2n-3+k$ edges.}, because
pointed pseudo-triangulations are {\it maximal} pointed sets edges
on any planar set of points and must have exactly $2n-3$ edges
(cf. \cite{streinu}). See Figure~\figref{pseudoT} for an example.



\section{Main result: Inductive Proof via Henneberg construction}
\label{section:main}

We are now ready to give our first proof of the Main Theorem, in
the following slightly more general form, and extend it to planar
Laman-plus-one graphs.

\begin{theorem}
\label{thm:mainL} Any plane Laman graph has a pointed
pseudo-triangular embedding.
\end{theorem}

\begin{theorem}
\label{thm:mainL1} Any plane Laman-plus-one graph has a
pointed-plus-one pseudo-triangular embedding.
\end{theorem}

Both proofs have the same structure and are divided into three
steps: topological, combinatorial and geometric. They are based on
the corresponding Henneberg constructions, the common theme of
this section. The easy topological lemmas shows that any Henneberg
construction on a (Laman or Laman-plus-one) plane graph can be
carried out in a manner compatible with the plane embedding. The
combinatorial lemmas construct a combinatorial
pseudo-triangulation while performing a topological Henneberg
construction. Finally, in the geometric lemmas we show how to
perform geometrically, with straight-line egdes, the Henneberg
extension (rather than just combinatorially). All the ideas are
already contained in the proof of the first theorem: the second
one will only be sketched, with indications of where the
differences lie.

\subsection{Pseudo-Triangular Embeddings of Plane Laman Graphs}

The proof of Theorem \ref{thm:mainL} is a consequence of the four
lemmas stated and proven below. Lemma \ref{lem:outerFace} reduces
the construction to the case when the outer face is a triangle.
Lemma \ref{lem:topoL} provides the framework for a Henneberg
induction on plane graphs. This is then used in Lemma
\ref{lem:cpt} to compute a combinatorial pseudo-triangulation
assignment and in Lemma \ref{lem:geom} to realize the same thing
geometrically. Theorem \ref{thm:mainL} follows from Lemma
\ref{lem:geom}. We remark that Lemma \ref{lem:cpt} is not needed
for the proof of Theorem \ref{thm:mainL}. It is however natural to
include it here because it makes use of the same Henneberg
technique (ubiquitous in this section). It also gives a better
intuition about the combinatorial structure of the many
possibilities involved in a complete proof by case analysis of
Lemma \ref{lem:geom}.

\begin{figure}[ht]
\begin{center}
\ \psfig{figure=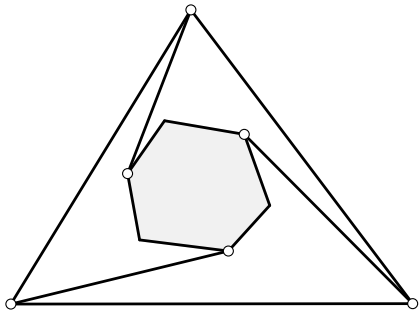,height=1.8in}
\end{center}
\caption{
Reducing to an embedding with a triangular outer face.}
\figlab{outerTriang}
\end{figure}

\begin{lemma} {\bf (Fixing the Outer Face)}
\label{lem:outerFace} Embedding a plane Laman graph as a
pseudo-triangulation reduces to the case when the outer face is a
triangle.
\end{lemma}

\begin{pf}
Let $G$ be a plane Laman  graph with an outer face having more
than three vertices. We construct another Laman  graph $G'$ of
$n+3$ vertices by adding $3$ vertices on the outer face and
connecting them into a triangle including the original graph. Then
we add an edge from each of the three new vertices to three
distinct vertices on the exterior face of $G$. See Figure
\figref{outerTriang}. We now realize $G'$ as a
pseudo-triangulation with the new triangle as the outer face. The
graph $G$, as a subgraph of $G'$, must be realized with {\it its}
outer face convex by the following argument. The three new
interior edges of $G'$ provide two corners each at their end-point
incident to the outer face and at least one corner in the interior
one. Since the three faces incident to them have nine corners in
total, the boundary of $G$ provides no corner to the three new
interior faces of $G'$.
\end{pf}

\smallskip

Notice that a planar Laman  graph (on $n$ vertices) always has at
least two triangular faces: the dual planar graph has $n-1$
vertices (including the vertex corresponding to the outer face)
and $2(n-1)-1$ edges, hence there exists at least two of degree
three. The construction in the previous lemma makes it possible to
use the stronger invariant of the Henneberg construction from
Lemma \ref{lemma:hennebergS} and to start any geometric embedding
with a triangular outer face, then to insert only on interior
faces. This feature is not needed in the proof of the topological
or combinatorial lemmas below.

\begin{figure}[ht]
\begin{center}
\ \psfig{figure=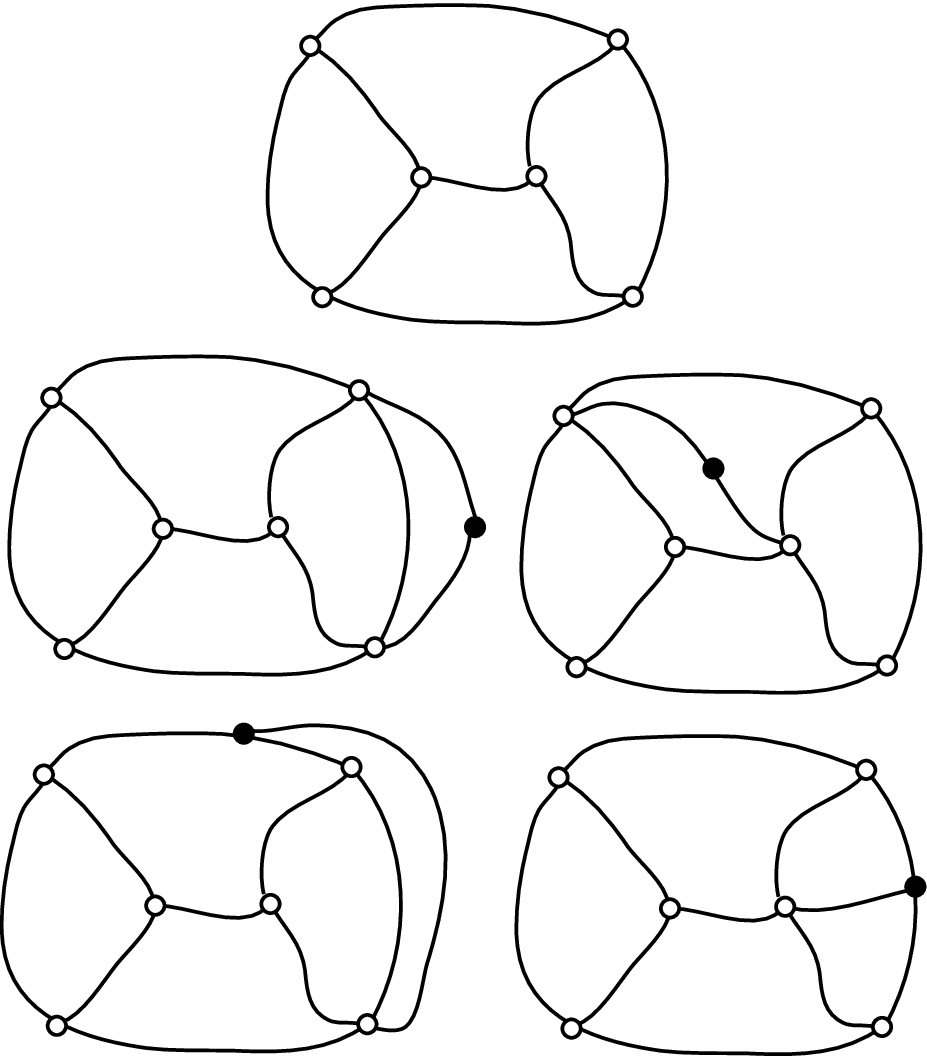,height=3.0in}%
\end{center}
\caption{
A plane Henneberg construction. Top row: $G_n$, to which a new
vertex will be added. Middle row: Henneberg I on the outer, resp.
an interior face. Bottom row: Henneberg II on the outer, resp. an
interior face.} \figlab{planarHenne}
\end{figure}

\begin{lemma} {\bf (The Topological Lemma)}
\label{lem:topoL}%
Every plane Laman graph has a {\em plane Henneberg construction}
in which:

\begin{enumerate}

\item All intermediate graphs are plane%
\item At each step, the topology is changed only on edges and
faces involved in the Henneberg step: either a new vertex is added
inside a face of the previous graph ({\em Henneberg I}), or inside
a face obtained by removing an edge between two faces of the
previous graph ({\em Henneberg II}).
\end{enumerate}

In addition, if the outer face of the plane graph is a triangle,
we may perform the Henneberg construction starting from that
triangle. The Henneberg steps will never insert vertices on the
outer face.
\end{lemma}

\begin{pf}  We follow the structure of the basic Henneberg construction from Lemma
\ref{lemma:hennebergS}. Find an appropriate vertex of degree $2$
or $3$. Removing it, and its incident edges, merges two (resp.
three) faces into one. The other endpoints of the removed edges
are incident to this face, hence the added edge in the Henneberg
II step simply splits this face and maintains the planarity of the
embedding.
\end{pf}

\begin{figure}[ht]
\begin{center}
\ \psfig{figure=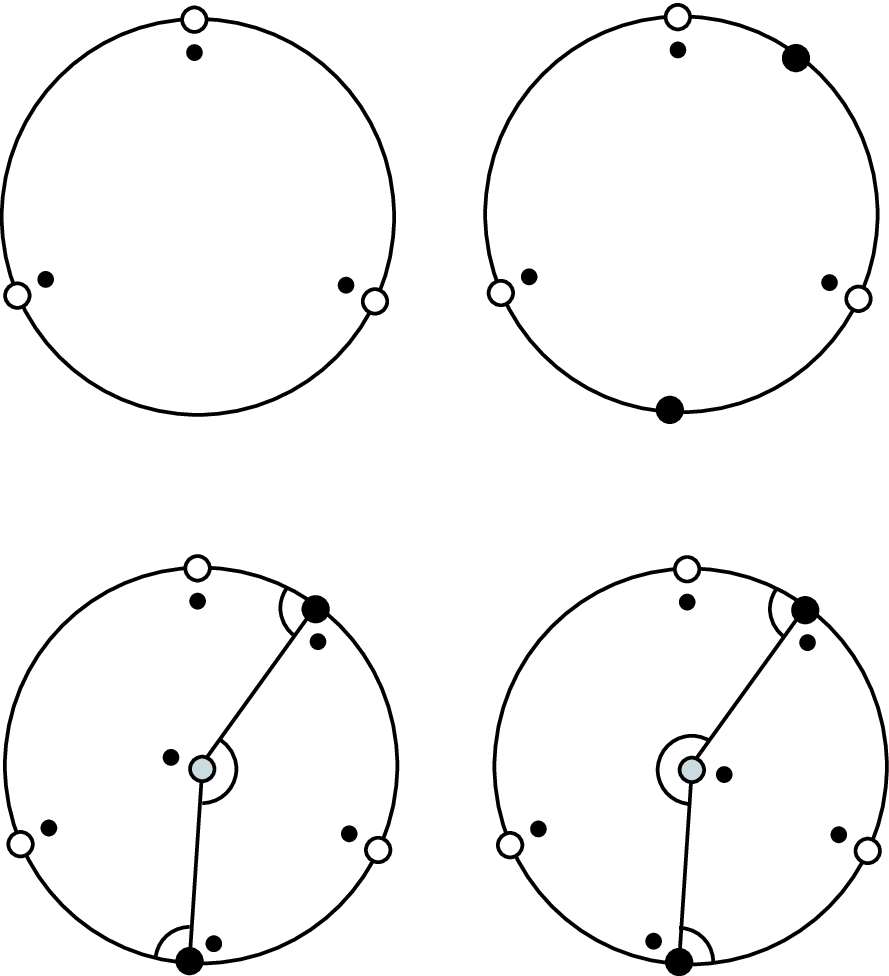,height=2.4in}
\end{center}
\caption{ Extending a combinatorial pseudo-triangulation in a
Henneberg I step. Top left: the combinatorial face is represented
as a circle with its three corners, denoted by white vertices,
marked {\it small} (a small black dot). Top right, a
representative situation: the two endpoints of the newly added
edges (in black) are distributed on two distinct side chains of
the face. Two distinct labelings are possible in this case (bottom
row): the newly created angles after the insertion of the new
vertex (grey) are labeled with a small black dot for {\em small}
(or {\em convex}) and with a large arc for {\em big} (or {\em
reflex}). } \label{fig:cpt1}
\end{figure}

\smallskip

\begin{lemma}{\bf (The Combinatorial Lemma)}
\label{lem:cpt} Every plane Laman graph admits a combinatorial
pseudo-triangulation assignment.
\end{lemma}

\begin{pf}
Let $G_n$ be a plane Laman graph on $n$ vertices. We may assume
that the outer face is a triangle. We proceed with a plane
Henneberg construction guaranteed by Lemma \ref{lem:topoL}, which
will insert only on interior\footnote{This is just a technical
simplification reducing the size of our case analysis. The reader
may verify that the Henneberg steps work as well for insertions on
the outer face.} faces. The base case is a triangle and has a
unique cpt labeling. At each step we have, by induction, a cpt
labeling which we want to extend. The proof will guarantee that
each one can be extended (so there is no need to backtrack).

\begin{figure}[ht]
\begin{center}
\ \psfig{figure=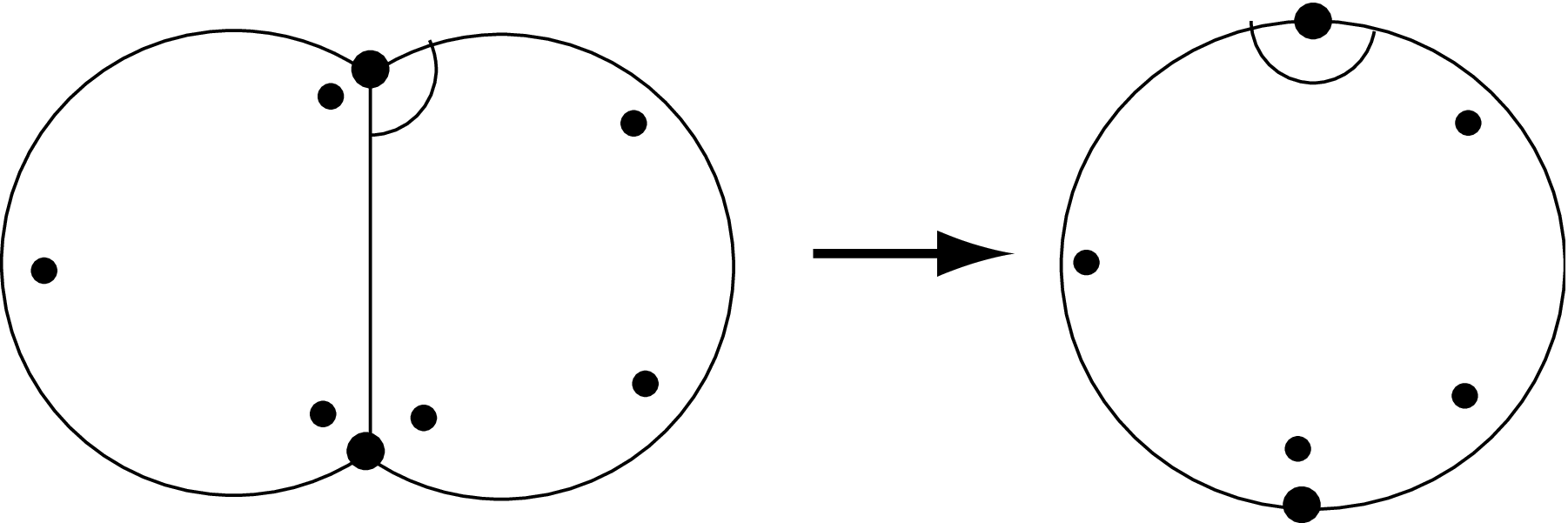,height=1.0in}
\end{center}
\caption{Merging two faces into one: a representative case for the
analysis of a combinatorial Henneberg II step. The markings of
{\em small} and {\em big} angles follow the conventions from
Figure \figref{cpt1}. The black vertices are the endpoints of the
removed edge $v_1v_2$.}%
\label{fig:cpt2A}
\end{figure}

In a Henneberg I step, the new vertex $v$ is inserted on a face
$T$ (already labeled as a pseudo-triangle), and joined to two old
vertices $v_1$ and $v_2$. The new edges $vv_1$ and $vv_2$
partition the face $F$ and its three corners into two. Three cases
may happen. Either the three corners fall one in one face and two
in another face, or one of them is split by a new edge (say $v_1$
is a corner) and the other two are either both in the same new
face or are separated, or two corners are split, and the third
corner is in one of the two newly created faces. In either case,
the assignment of {\it big} and {\it small} labels is what one
would expect: a small angle is split into two small angles, a big
angle is split into a big and a small angle, and the new point
gets exactly one big angle. We illustrate one representative case
in Figure \ref{fig:cpt1} and leave the rest of the straightforward
details of this case analysis to the reader.

\begin{figure}[ht]
\begin{center}
\ \psfig{figure=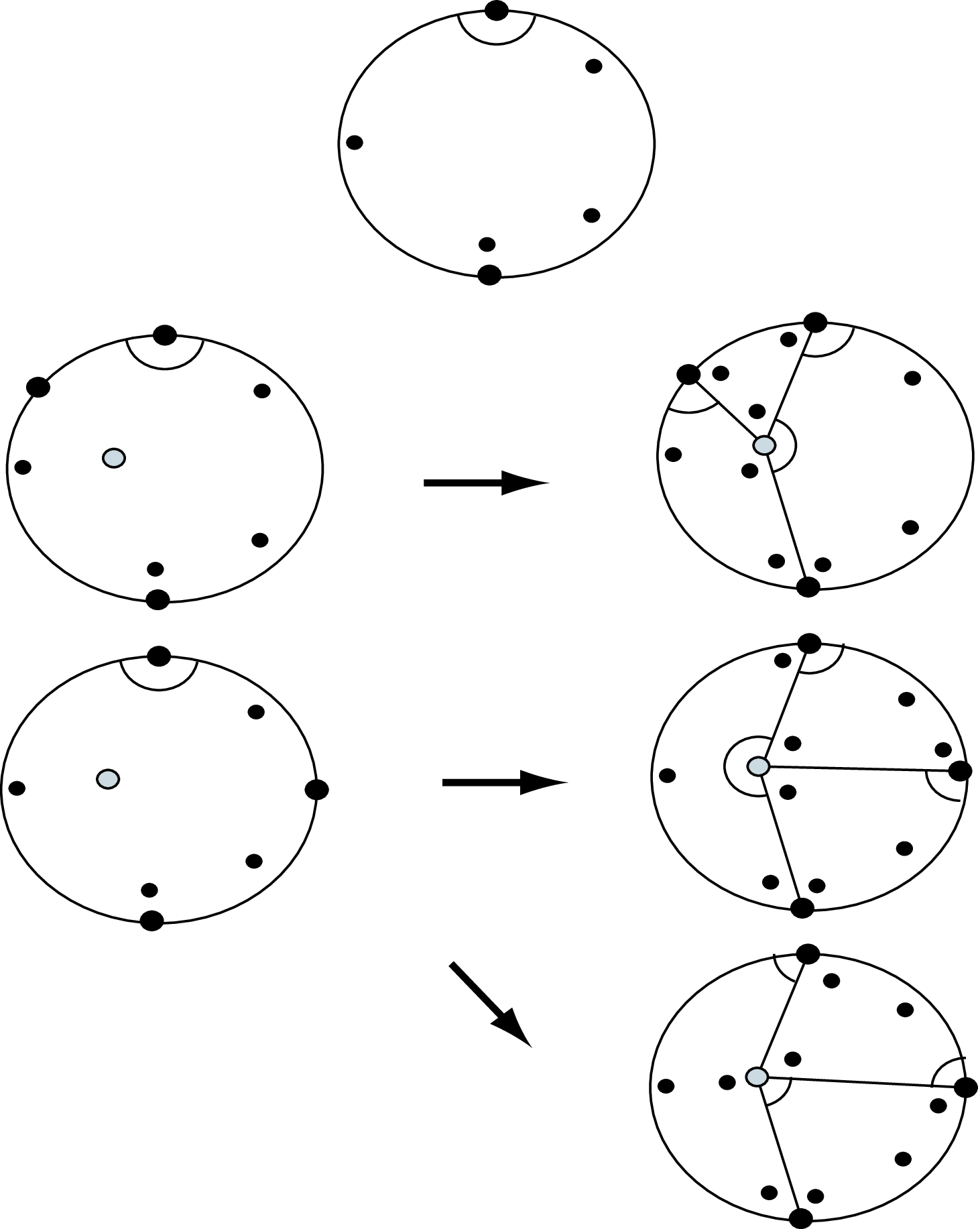,height=3.2in}
\end{center}
\caption{ Extending a combinatorial pseudo-triangulation in a
Henneberg II step. Top: a representative case of a combinatorial
face with four corners. Left: two (out of three) possible
placements of the third vertex. Right: the possible
labelings of the induced faces as combinatorial pseudo-triangles.}%
\label{fig:cpt2B}
\end{figure}

In a Henneberg II step, an edge $v_1v_2$ is first removed, merging
two faces labeled as combinatorial pseudo-triangles into one face
$T$. In this process, some angles are merged into one: their
labels must be reassigned, but we make no changes to the labels of
the other angles. The rules for assigning labels to merged angles
are simple, mimicking what one would expect to happen in a
straight-line situation: if one old angle was {\it big}, the
merged angle is marked {\it big}, otherwise {\it small}. The face
$T$ thus gets exactly four small angles. Its boundary is separated
by the vertices $v_1$ and $v_2$ into two chains: each contains at
least one corner. Four cases may happen: $v_1$ and $v_2$ are both
{\it small} (corners), separating the other two corners; only one
is {\it small} (say, $v_1$), and the other corners are distributed
as $1$-$2$ on the chains; or both $v_1$ and $v_2$ are {\it big},
and the four corners are distributed as either $2$-$2$ or $1$-$3$.
See Figure \figref{cpt2A} for a representative case (the other are
similar and are left to the reader). Notice that it is impossible
to have all four corners on only one chain induced by $v_1$ and
$v_2$.

The new vertex $v$ is now inserted inside this face $T$, and
joined to the old vertices $v_1$ and $v_2$, and to some other
vertex $v_3$ on $T$. The new edges $vv_1$, $vv_2$ and $vv_3$
partition the face $T$ and its four corners into three parts,
which can be assigned the labels in several ways. See Figure
\figref{cpt2B} for a representative case: the systematic
verification  of all the cases is straightforward and is left to
the reader.
\end{pf}


\medskip

Notice that in general the pointed combinatorial
pseudo-triangulation guaranteed by Lemma \ref{lem:cpt} is not
unique. The Lemma shows in fact how to systematically generate all
of them in a non-backtracking manner: if we succeeded in finding a
pointed cpt at step $n$, we simply extend it at the next Henneberg
step. We next prove that at least one of them is realizable with
straight-lines via a similar Henneberg extension technique.

\smallskip

\begin{lemma}{\bf (The Geometric Lemma)}
\label{lem:geom} Every plane Laman graph $G$ can be embedded as a
pseudo-triangulation.
\end{lemma}

\begin{pf} Let $G_n$ be a plane Laman graph on $n$ vertices with a
triangular outer face. Assume we have a plane Henneberg
construction for $G_n$ starting with the outer face and adding
vertices only on interior faces. We basically follow the same
analysis as in Lemma \ref{lem:cpt}. This time, however, we will
not {\em choose} the big/small labels of the angles, but rather
show that {\em there exists} a way of placing a point $p_n$ inside
a face which realizes a compatible partitioning of the face into
pseudo-triangles as prescribed by the Henneberg step on the vertex
$v_n$ of degree $2$ or $3$.

As in Lemma \ref{lem:cpt}, the Henneberg I step is straightforward
on an interior face (which is what we do here). As an exercise
pointing out to the difference between the combinatorial and the
geometric case, we leave it to the reader to verify that this is
not the case on an outer face, where the placement of a vertex at
step $n+1$ may be constrained by the realization up to step $n$,
and thus may not directly allow an embedding with a certain
prescribed outer face.

The analysis of a Henneberg II step is identical to that performed
in the combinatorial lemma, and leads to several cases to be
considered. We illustrate here only a representative case (but
have verified them all). The important fact is that {\em it is
always possible} to realize with straight-lines {\em at least one}
of the possible Henneberg II combinatorial pseudo-triangular
extensions.

Consider the (embedded) interior face $F$ with four corners
obtained by removing an interior edge $p_ip_j$, and  let $p_k$ be
a vertex on the boundary of $F$. We must show that there exists a
point $p$ inside $F$ which, when connected to $p_i$, $p_j$ and
$p_k$ partitions it into three pseudo-triangles and is itself
pointed. The three line segments $pp_i$, $pp_j$ and $pp_k$ must be
tangent to the side chains of $F$. We define the {\it feasibility
region} of an arbitrary point $p_a$ on the boundary of $F$ as the
(single or double) wedge-like region inside $F$ from where
tangents to the boundary of $F$ at $p_a$ can be taken. The
feasibility region of several points is the intersection of their
feasibility regions. An important fact is that the feasibility
region of $p_i$ and $p_j$ always contains the part of the
supporting line of the removed edge $p_ip_j$, and that the
feasibility region of {\it any} other vertex $p_k$ cuts an open
segment on it. In fact, the feasibility region of $p_k$ intersects
the feasibility region of $p_i$ and $p_j$ in a non-empty feasible
$2$-dimensional region on one side or the other (or both) of this
segment. One can easily see that not only is this region
non-empty, but it contains a subregion where a placement of $p$
{\em as a pointed vertex} is possible (we call that a {\em
pointed-feasible region}). We skip the rest of the straightforward
details. See Figure \figref{geom} for a representative case.
\end{pf}

\begin{figure}[ht]
\begin{center}
\ \psfig{figure=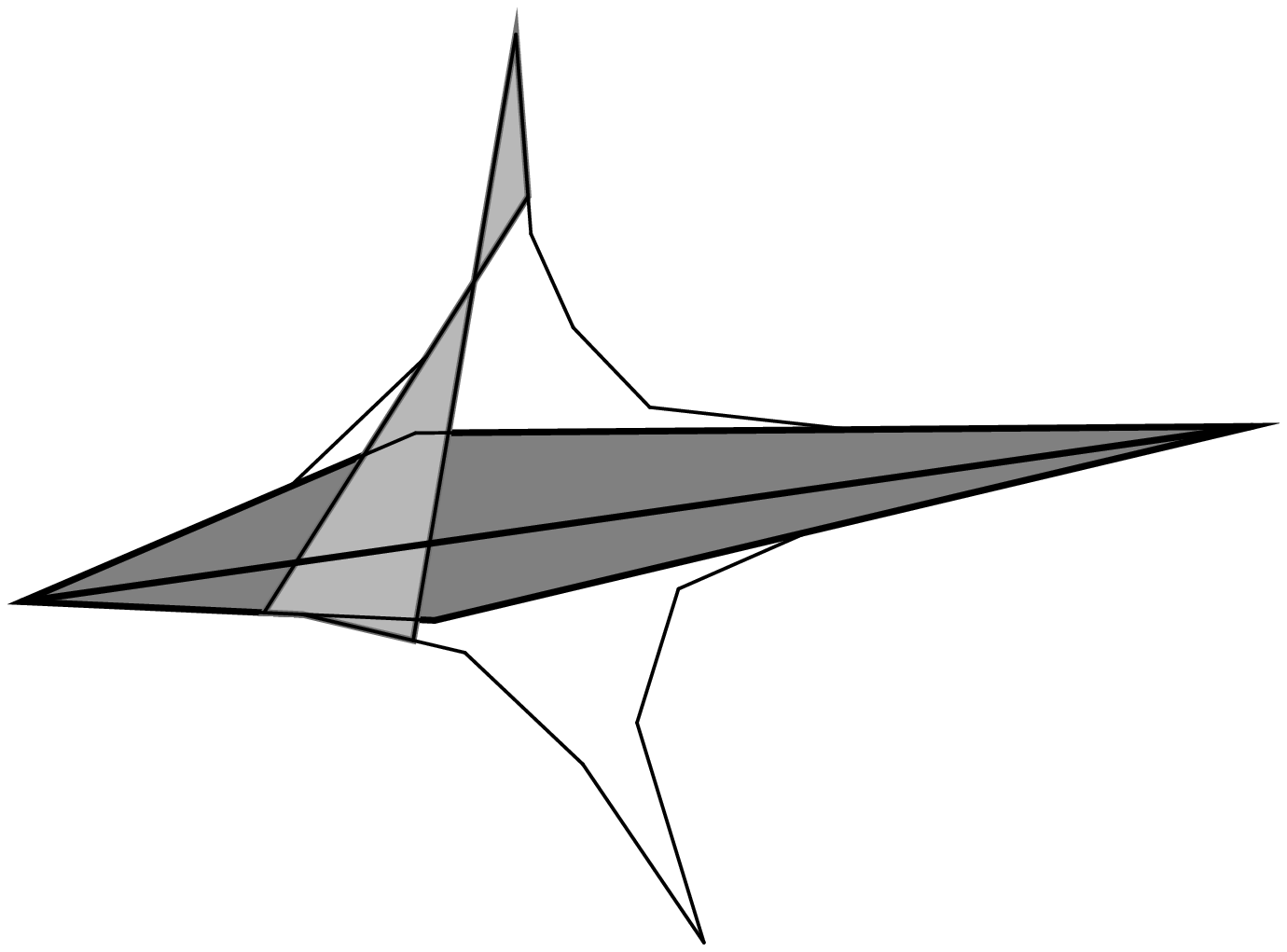,width=1.8in}
\psfig{figure=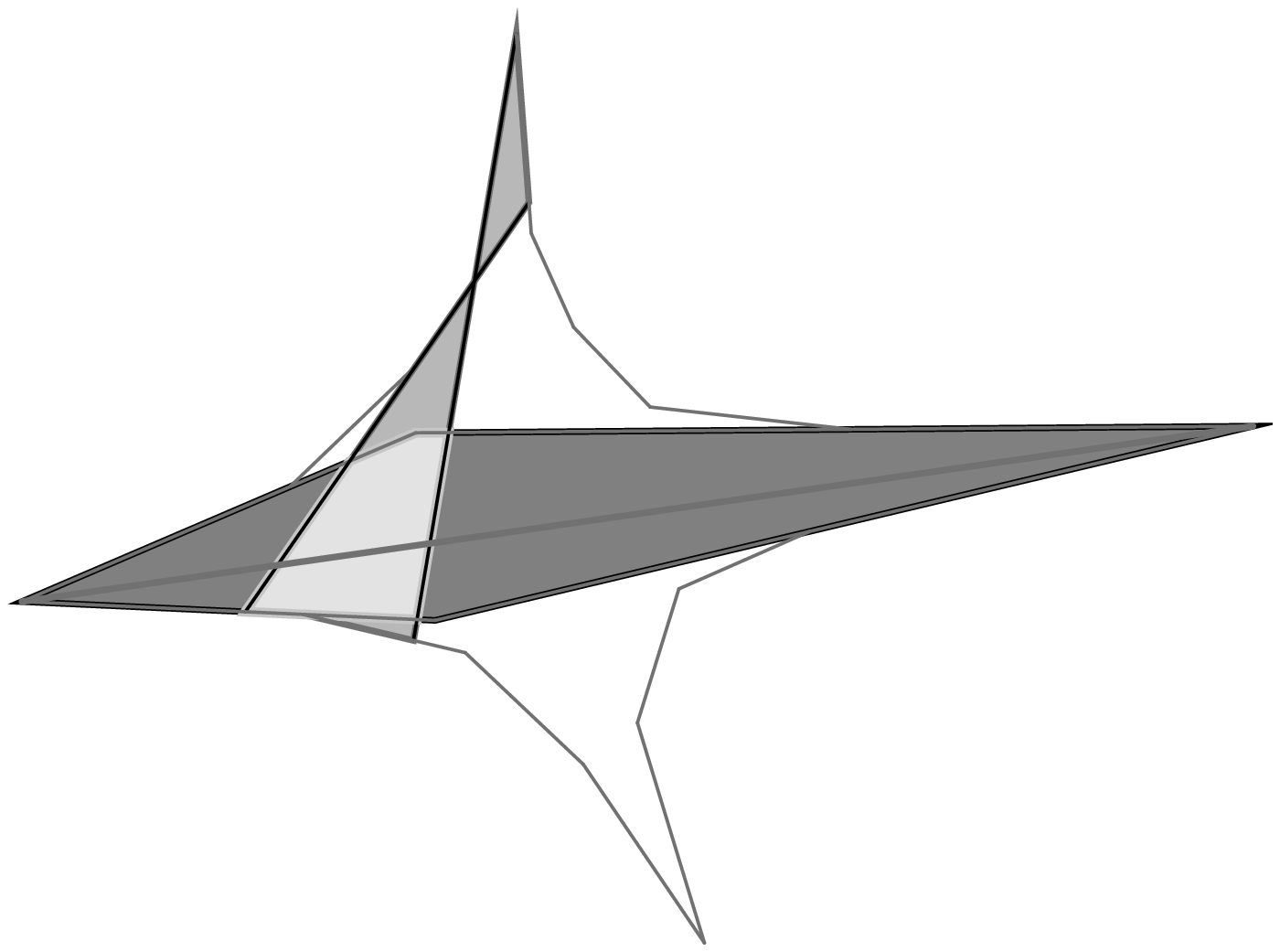,width=1.8in}
\ \psfig{figure=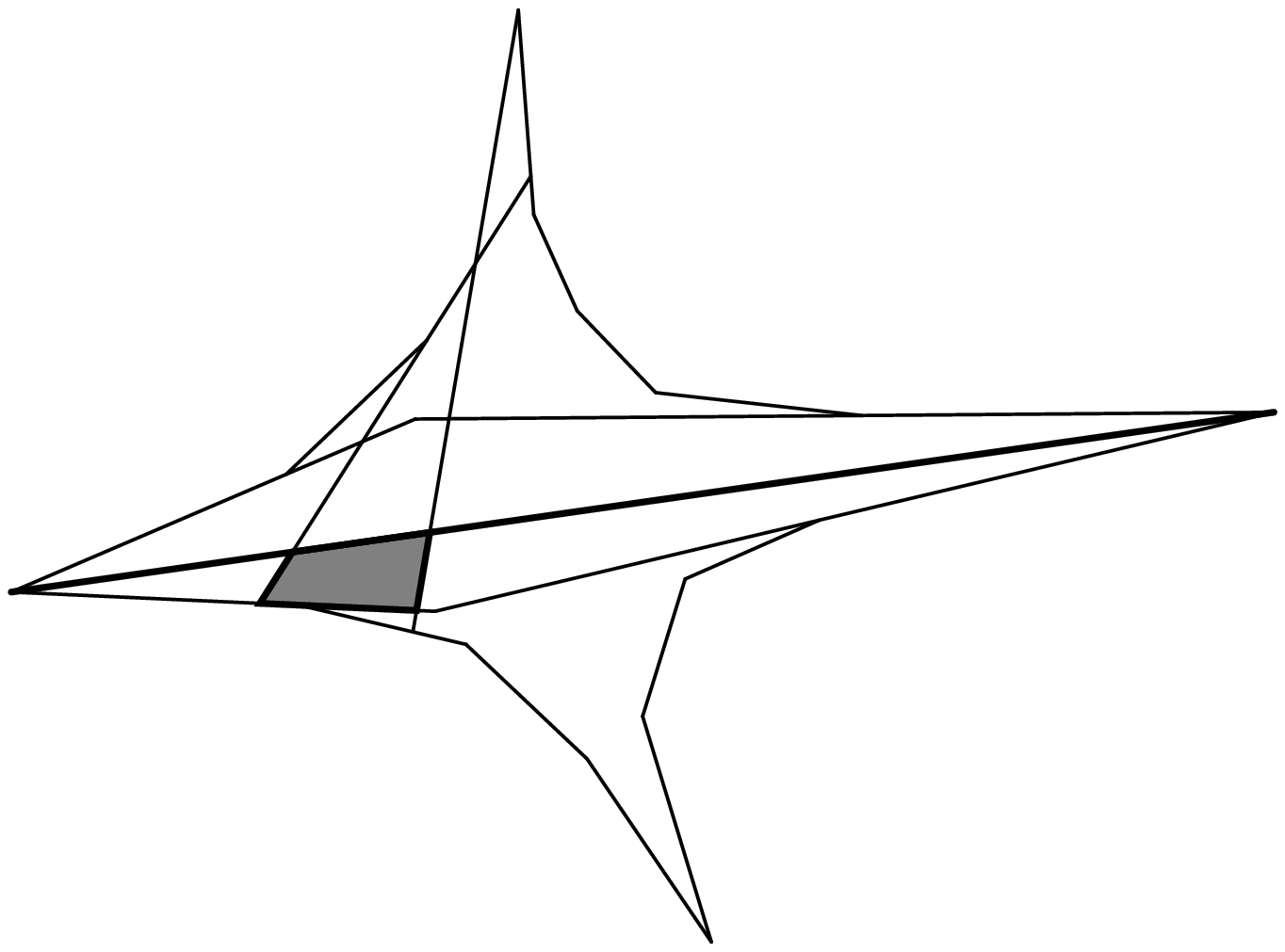,width=1.8in}
\psfig{figure=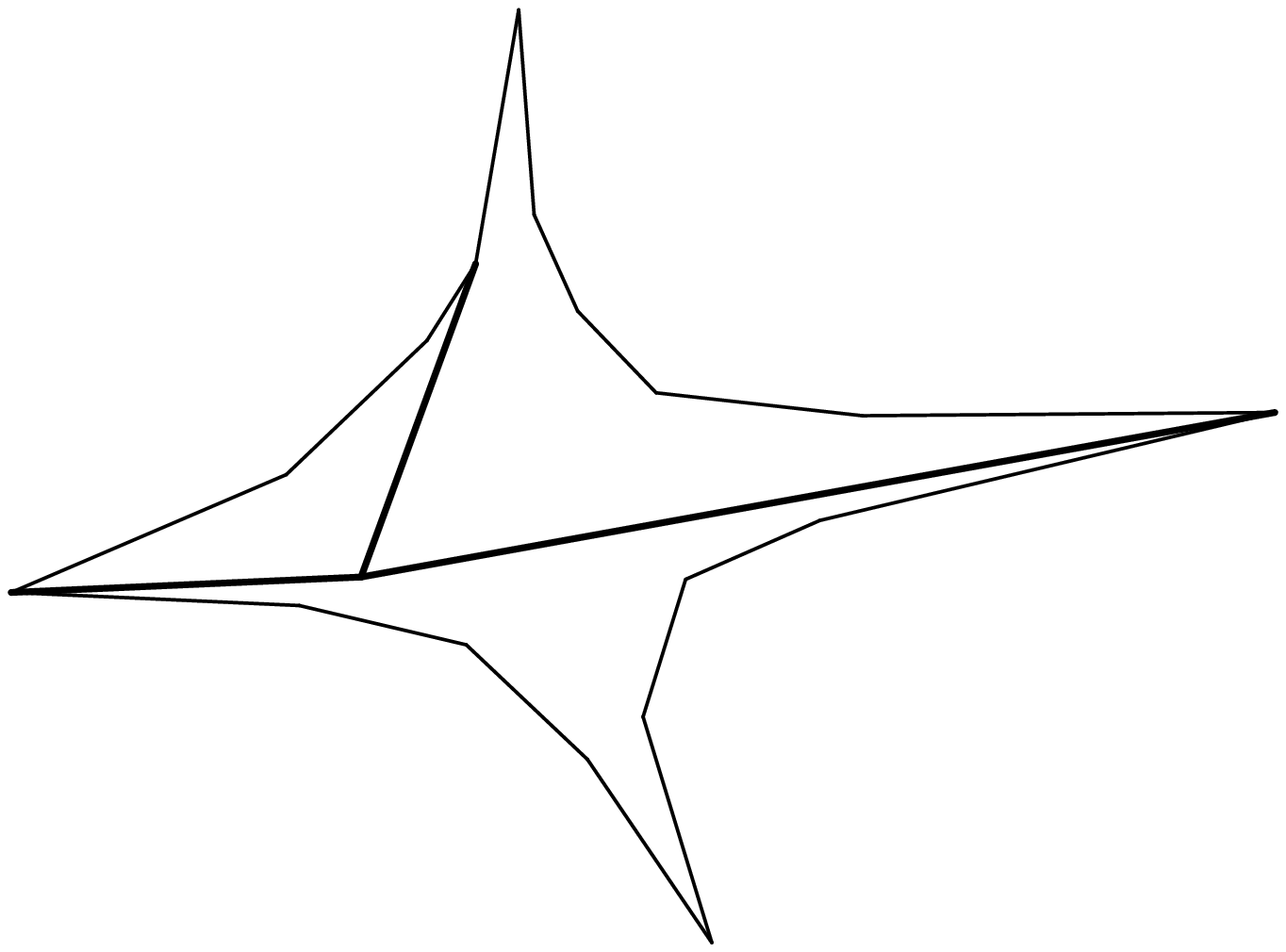,width=1.8in}
\end{center}
\caption{
Henneberg II step on an interior face. One sees that the feasible
region of $p_k$ (light grey) intersects the feasible region of the
two endpoints $p_i$ and $p_j$ (dark grey) of the removed edge. We
show the two feasible regions (of the pair $p_i, p_j$, resp.
$p_k$), their intersection, the final pointed-feasible region and
a placement of
a pointed vertex $p$ and its three tangents.}%
\figlab{geom}
\end{figure}

\medskip
\begin{proofT} \ref{thm:mainL}.
Let $G$ be a plane Laman graph. If its outer face is not a
triangle, apply Lemma \ref{lem:outerFace} to get a new graph $G_n$
which will contain $G$ in its embedding. Follow the plane
geometric Henneberg construction described in Lemma \ref{lem:geom}
to embed $G_n$ starting from a triangle and always inserting new
vertices in some interior face.
\end{proofT}


\medskip
\noindent {\bf Algorithmic analysis.} The proof of Theorem
\ref{thm:mainL} can be turned into an efficient polynomial time
algorithm. Given a Laman graph, verifying its planarity and
producing a plane embedding (stored as a quad-edge data structure
\cite{guibasstolfi} with face information) can be done in linear
time \cite{ht}. One then chooses an outer face and in linear time
one can perform the construction from Lemma \ref{lem:outerFace} to
get a triangular outer face. For producing a topological Henneberg
construction, we'll keep an additional field in the vertex data
structure, storing the degree of the vertex. We will keep the
vertices in a min-heap on the degree field. To work out the
Henneberg steps in reverse we need to do efficiently the following
operations: a) detect a vertex of minimum degree (which will be
$2$ or $3$), b) if the minimum degree is $3$, corresponding to a
vertex $v$, we must find an edge $e$ that will be put back in
after the removal of the neighbors of $v$, and c) restore the
quad-edge data structure. Step a) can be done in $O(\log n)$ time.
Step c) can be done in constant time. Step b) requires deciding
which of the three possibilities for $e$ among $v_1v_2$, $v_1v_3$
and $v_2v_3$ (where $v_1, v_2$ and $v_3$ are the neighbors of the
vertex removed in a reverse Henneberg II step) produces a Laman
graph. Testing for the Laman condition on a graph with $2n-3$
edges can be done by several algorithms (the algorithms of Imai
and Sugihara, via reductions to network flow or bipartite
matching, or via matroid (tree) decompositions, see \cite{walterS}
and the references given there), and takes $O(n^2)$. Therefore the
time for performing a reverse Henneberg step is dominated by b),
which gives a total running time of $O(n^3)$.

The embedding is done now by performing the Henneberg steps,
starting with the outer triangular face embedded on an arbitrary
initial triple of points. It is straightforward to see that each
step takes constant time to determine a position for the new
vertex, and the whole embedding takes linear time once the
Henneberg sequence is known.

The time complexity would be improved by a positive answer to the
open questions \ref{open:LamanPlanar} and \ref{open:edgeBack}
listed in the concluding Section \ref{section:open}.


\subsection{Pseudo-Triangular Embeddings of Plane Laman-plus-one Graphs}
\label{subsection:circuit}

We now turn to a proof of Theorem \ref{thm:mainL1} using Henneberg
constructions for plane Laman-plus-one graphs. It is very similar
to the proof of Theorem \ref{thm:mainL}. It is instructive though
to see the differences, which lie in the combinatorial (and hence
also geometric) pseudo-triangulation assignment, where we must
keep track of the non-pointed vertex. We have two items which may
in principle be prescribed: the outer face and the vertex to
become the unique non-pointed one. The non-pointed vertex may only
be interior to the circuit. In a Henneberg construction, we will
see that it is easy to prescribe either the outer face or the
interior vertex to be non-pointed, but the analysis becomes more
complicated for the prescription of both. In Section
\ref{section:cpt} we use a different, global argument to do the
simultaneous prescription of the outer face and of the non-pointed
vertex, in the case of a rigidity circuit.

The next two lemmas are straightforward extensions of the Laman
case.

\begin{lemma} {\bf (Fixing the Outer Face)}
\label{lem:outerFace1} Embedding a plane  Laman-plus-one graph as
a pseudo-triangulation reduces to the case when the outer face is
a triangle.
\end{lemma}

Notice that this would fail if we only planned to deal with a
circuit.

\begin{lemma} {\bf (The Topological Lemma)}
\label{lem:topoL1}%
Every plane  Laman-plus-one graph has a {\em plane Henneberg
construction}.
\end{lemma}


A planar Laman-plus-one  graph (on $n$ vertices) always has at
least two triangular faces: the dual planar graph has $n$ vertices
(including the vertex corresponding to the outer face) and $2n-2$
edges, hence there exists at least two of degree three. The
previous construction allows to prescribe the outer face in a
geometric embedding, should we want to do so, and is not needed in
the proof of the  combinatorial lemma below.

\begin{lemma}{\bf (The Combinatorial Lemma)}
\label{lem:cpt1} Every plane Laman-plus-one graph admits a
pointed-plus-one combinatorial pseudo-triangulation assignment.
\end{lemma}

\begin{pf}
The proof has the same basic structure (but more cases to analyze)
as Lemma \ref{lem:cpt}, and relies on the details of the Henneberg
construction from Lemma \ref{lemma:hennebergC}. The base case is
$K_4$ which has a unique cpt assignment for a choice of an outer
face. It is easy to see that Henneberg I steps cause no problem,
and the Henneberg II steps work as before when the vertex of
degree $3$ is not inside the circuit, is not the pointed vertex,
and it is not incident to it.

Let $v_iv_j$ be the removed edge and $v_k$ the third vertex
involved in the Henneberg II step. The only problematic case is
when the edge $v_iv_j$ is incident to the unique non-pointed
vertex. In this case, the resulting face after the removal of
$v_iv_j$ is a pseudo-triangle: it has three, not four corners. We
must argue that {\em in at least one combinatorial
pseudo-triangulation} compatible with the information so far, the
three vertices $v_i$, $v_j$ and $v_k$ cannot lie all three on the
same side-chain of this face, otherwise the extension to a cpt is
impossible.

There is a way around this, which would guarantee that both the
outer face and the non-pointed vertex could be prescribed. We will
describe this, in a more general setting, in a forthcoming paper.
For the time being, it suffices to notice that if this happens
during the Henneberg construction, we can always pick up one of
the other three vertices guaranteed to have degree three (when
there are no degree two vertices), and continue from there. Notice
that this may change the outer face assignment, though.
\end{pf}

\begin{lemma}{\bf (The Geometric Lemma)}
\label{lem:geom1} Every plane Laman-plus-one graph $G$ can be
embedded as a pseudo-triangulation.
\end{lemma}

This proof, and the proof of Theorem \ref{thm:mainL1} are now
straightforward extensions of those done for the Laman case.
Notice that it may not be possible in general to guarantee a
certain outer face or non-pointed vertex.

\medskip
\noindent {\bf Remarks.} The inductive technique described in this
section works in fact for {\it rigid} graphs on $n$ vertices and
fewer than $2n$ edges (i.e. Laman graphs with at most two extra
edges). Indeed, the only ingredient that is needed is the
existence of a vertex of degree at most $3$ whose removal either
leaves a rigid graph or a graph with one degree of freedom. Given
a rigid graph with a combinatorial pseudo-triangulation assignment
we can add edges and combinatorially assign the big and small
labels to the angles, preserving at each step the property that
every face has exactly three small angles and every vertex has at
most one big angle and the outside face has only big angles.
However, not all of these combinatorial angle assignments are
geometrically realizable.

\smallskip
The {\bf algorithmic analysis} is similar to the case of plane
Laman graphs.


\section{Combinatorial pseudo-triangulations}
\label{section:cpt}

In this section we extend the results from Section
\ref{section:main} on pointed and pointed-plus-one combinatorial
pseudo-triangulations. We present a global, non-inductive
technique for generating cpt assignments for planar Laman graphs
and planar circuits. It is based on a reduction to finding perfect
matchings in a certain associated bipartite graph. By showing that
Hall's condition is satisfied, we are guaranteed to have a
solution (and hence a cpt) for both plane Laman graphs and
circuits. We also show that the existence of a pointed
combinatorial pseudo-triangulation assignment is not restricted to
plane Laman graphs or circuits.

\smallskip

Let $G=(V,E,F)$ be a plane graph with vertices $V$, edges $E$ and
faces $F$. Assume $|V|=n$ and $|E|=2n-3$. Euler's relation implies
that $|F|=n-1$. Denote by $F'$ the set of interior faces and by
$f_o$ the outer face (with $h$ vertices), $F=F'\cup\{f_o\}$. We
define a bipartite graph $H$ with the two sets of the bipartition
labeled $V$ and $W$. $V$ stands for the set of vertices $V$ of $G$
and has $n$ elements. The set $W$ corresponds to the faces $F$ of
$G$ taken with certain multiplicities. For an interior face $f\in
F'$ of degree (number of edges on the face) $d_f$, we will put
$d_f-3$ vertices in $W$. For the outer face $f_o$ we will put
$h=d_{f_o}$ nodes in $W$. The total number of elements in $W$ is
thus $\sum_{f\in F'}(d_f-3) + h = \sum_{f\in F} d_f - 3 |F'|=2 |E|
- 3(n-2) = 2(2n-3) - 3(n-2) = n$.

\begin{figure}[ht]
\begin{center}
\ \psfig{figure=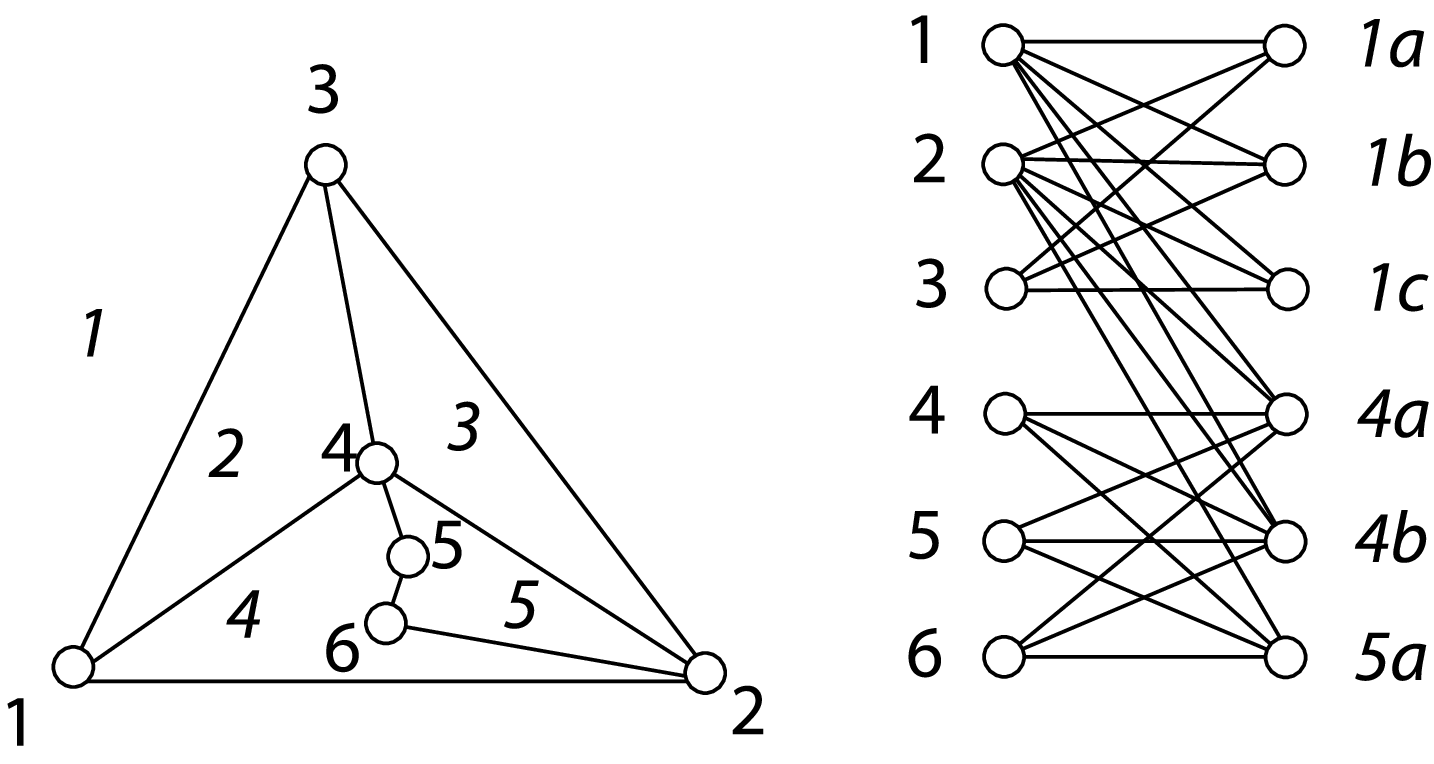,width=3.2in}
\end{center}
\caption{A plane graph with $2n-3$ edges and its associated
bipartite graph $H$.}%
\figlab{match}
\end{figure}

A vertex $v\in V$ is connected in $H$ to the vertices in $W$
corresponding to the interior faces $f$ of degree larger than $3$
to which it belongs in $G$, and to the vertices corresponding to
the outer face (if it belongs to it). Hence if $v$ belongs to
three faces $f_1, f_2, \cdots $, and these faces have
multiplicities $d_1, d_2, \cdots$ in $W$, then $v$ is connected to
$d_1$ copies of the vertex for $f_1$, $d_2$ copies for $f_2$, etc.
See Figure \figref{match}. The $6$ vertices and $5$ faces of
degrees $3$ (outer face $1$), $3$ (faces $2$ and $3$), $4$ (face
$5$) and $5$ (face $4$) lead to the bipartition sets
$V=\{1,2,3,4,5,6\}$ and $W=\{1a, 1b, 1c, 4a, 4b, 5a \}$, connected
by edges as in the figure.

The connections (edges) in the bipartite graph $H$ represent
potential assignments of {\it big} angles, where an {\it angle} is
viewed as a pair {\it (vertex, face) }. Since each vertex must
receive a big angle, we want a perfect matching. Since each
interior face receives all but three big angles, and the outer
face receives all big angles, the choice of multiplicities
reflects just that. These considerations lead to the following
Lemma.

\begin{lemma}
\label{lem:match}%
There is a one-to-one correspondence between the combinatorial
pseudo-triangulations of a plane graph $G$ with $n$ vertices and
$2n-3$ edges and the perfect matchings in the associated bipartite
graph $H$.
\end{lemma}

\begin{figure}[ht]
\begin{center}
\ \psfig{figure=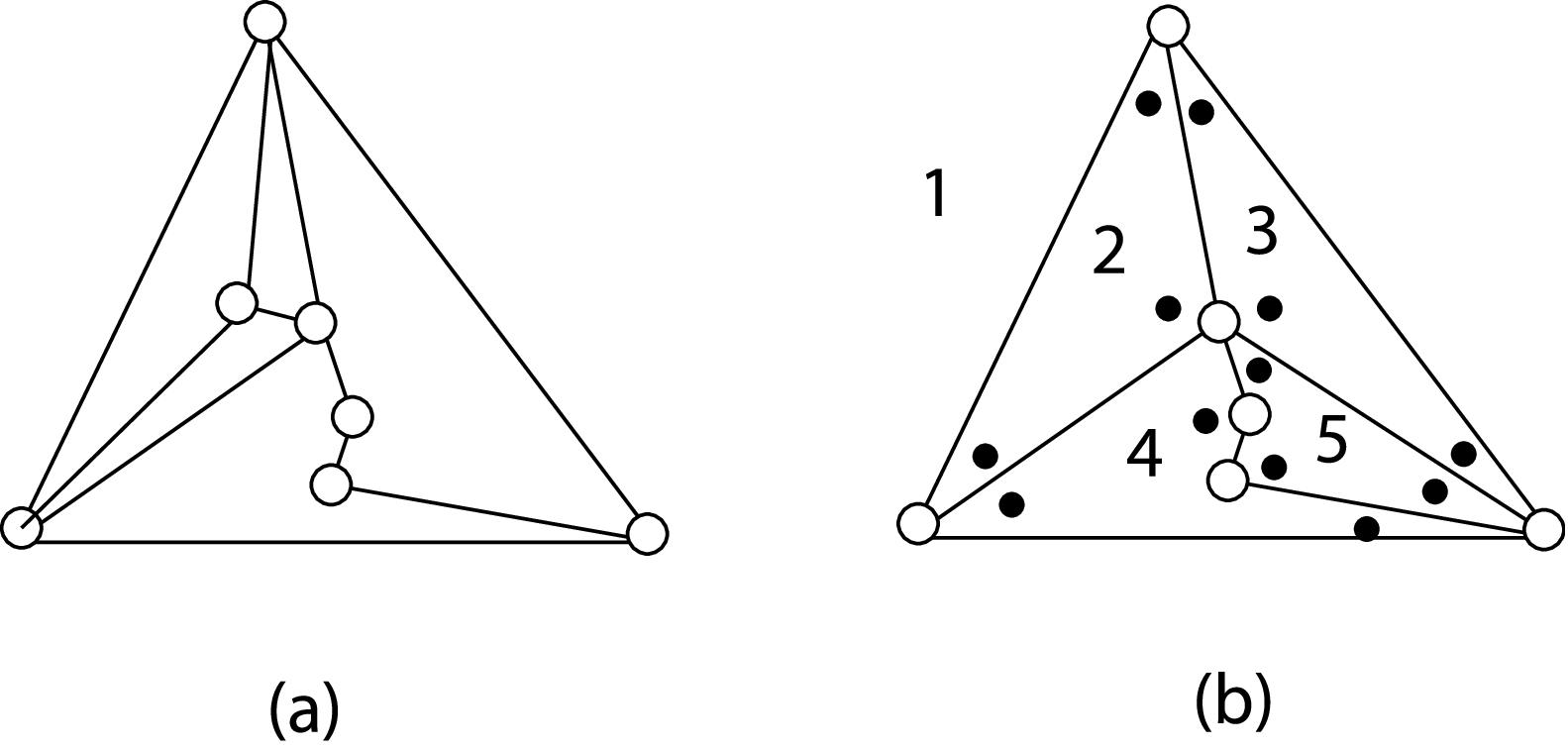,width=3.4in}
\end{center}
\caption{
(a) A plane graph with $2n-3$ edges and no combinatorial pseudo-triangulation assignment.
(b) A plane non-Laman graph with a cpt assignment.}%
\figlab{cptN}
\end{figure}

In general, plane graphs satisfying the conditions of the previous
Lemma may or may not have combinatorial pseudo-triangulation
assignments. See Figure \figref{cptN} for examples. But for Laman
graphs, we are guaranteed a solution. The main  result of this
section is:

\begin{theorem}
\label{thm:match}%
If $G$ is a Laman graph, then $H$ has a perfect matching. Hence
$G$ has a pointed combinatorial pseudo-triangulation.
\end{theorem}

\begin{pf}
We will check Hall's condition to guarantee the existence of a
perfect matching. Let $A\subset V$ be a subset of $|A|=a$
vertices. Let $F_A$ of size $|F_A|=f_a$ be the set of faces
incident to the vertices in $A$, and let $D=\sum_{f\in F_A}d_f$.
We need to show that $a\leq D-3 f_A$.

It suffices to carry out the analysis on different face-connected
components of $F_A$ separately. See Figure \figref{faces}. One
face-connected component is a polygon with (say) $b$ boundary
edges, $b' \leq b$ boundary vertices, $h\geq 0$ holes and $e_i$
interior edges. We have $D=2 e_i +b$. By Euler's relation $(a+b')
+ (f_A + h + 1) = (b + e_i) + 2$. Hence $f_A = e_i - a + 1 + b -
b' - h$, where $\Delta := b - b' - h \geq 0$. Laman's condition
implies that $b + e_i \leq 2(a+b') - 3$, hence $e_i \leq 2a +
2b'-b-3$. Now to show $a+3f_A \leq D$ we need $a + 3(e_i - a + 1 +
b - b' -h) \leq 2 e_i + b$, i.w. $e_i \leq 2 a - 3 + 3 b' - 2b +
3h$. Since we know $e_i \leq 2a + 2b'-b-3$, it remains to show
$2b' - b \leq 3b' - 2b + 3h$, i.e. $b \leq b'+3h$, which is
obviously true.
\end{pf}

\begin{figure}[ht]
\begin{center}
\ \psfig{figure=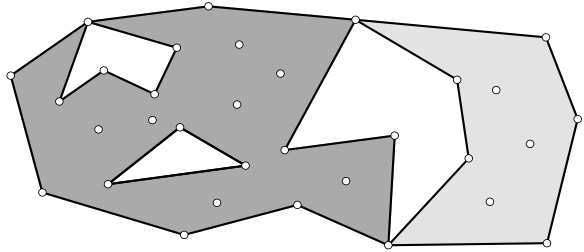,width=3.2in}
\end{center}
\caption{The analysis in the proof of Theorem \ref{thm:match}: it
suffices to analyze separately the face-connected components of $F_A$
(two in this case, shaded slightly differently).}%
\figlab{faces}
\end{figure}

The result of Theorem \ref{thm:match} extends to the case of plane
circuits. Moreover, we will be able to show in this case a more
general version of Theorem \ref{lem:topoL1}, by being able to
prescribe both the outer face and the non-pointed vertex (which
may be chosen as any vertex non-incident to the outer face). In
this case, the associated bipartite graph is slightly different:
the set $V$ contains all vertices but one, namely the vertex
prescribed to be the non-pointed one. The set $W$ has the same
description, but its size now is $n-1$ (because the number of
faces is $n$).

\begin{theorem}
\label{thm:matchC}%
If $G$ is a plane circuit, then $H$ has a perfect matching. Hence
$G$ has a pointed-plus-one combinatorial pseudo-triangulation with
a prescribed outer face and prescribed non-pointed vertex.
\end{theorem}

\begin{pf}
The analysis from proof of Theorem \ref{thm:match} still holds,
because in the case of a circuit, the condition on subsets of size
$k<n$ is {\it exactly} the same as for Laman graphs: they span at
most $2k-3$ edges. Therefore the analysis works whenever $F_A$
does not cover the whole polygon. Since at least one vertex is
missing from $A$ (the vertex prescribed to be non-pointed), this
is always the case.
\end{pf}

\medskip
\noindent {\bf Algorithmic analysis.} To check whether a graph
admits a combinatorial pseudo-triangulation (and to compute one)
we will use the $O(n^{3/2})$ time algorithm for the maximum flow
problem of Dinits (see~\cite{t-dsna-83}) to solve the bipartite
matching problem described above.

\smallskip
\noindent {\bf Remark.} This set of degree-constrained subgraphs
of a bipartite graph can be modelled as a network flow problem.
Thus the set of combinatorial pseudo-triangulations of a given
graph (with a given planar embedding, including a specification of
the outer face) is in one-to-one correspondence with the vertices
of a polytope, given by the equations and inequalities of the
network flow.


\section{Stretching combinatorial pseudo-triangulations}
\label{section:proofT}

We have seen in the previous sections
two proofs of the fact that every plane Laman graph can be
assigned a combinatorial pseudo-triangulation labeling. The
technique from Section~\ref{section:main} does not realize
geometrically every such possible combinatorial structure. In this
section we give the strongest version of the main result by
proving the following theorem.

\begin{theorem}
\label{thm:mainC} For any plane {\em Laman graph} $G$ and for any
of its combinatorial pseudo-triangular assignments, there exists
(and can be efficiently found) a compatible straight-line
embedding. The same holds for plane {\em circuit graphs}.
\end{theorem}

The proof is a consequence of two general results of independent
interest. We first give in Theorem \ref{thm:stretchable} two
characterizations of stretchable combinatorial
pseudo-triangulations. The stretchability proof relies on a
directed version of Tutte's Barycentric Embedding Theorem,
(Theorem \ref{thm:directedTutte}). Finally, we show that the
characterization in Theorem \ref{thm:stretchable} is satisfied for
Laman (Theorem \ref{thm:Lamanstretch}) and circuit plane graphs
(Theorem \ref{thm:circuitstretch}) with cpt assignments.


\subsection{Two characterizations of stretchability}

In this section we give two combinatorial characterizations of
stretchability of combinatorial pseudo-triangulations in terms of
the number of {\it corners} of planar subcomplexes and in terms of
$3$-connectivity properties of an associated directed graph.

\smallskip

Let $G=(V,E)$ be a plane graph with a combinatorial
pseudo-triangulation labeling. We do not impose any restrictions
on its number of non-pointed vertices or rigidity properties. As a
plane graph, every subgraph $G_S=(S,E_S)$ induced by a subset of
vertices $S\subset V$ has an induced plane embedding and a
well-defined unbounded region. The {\it boundary} of the unbounded
region
consists of cycles of vertices and edges, with one cycle for each
connected component of $G_S$. Some edges and/or vertices may be
repeated in these cycles. For example, if $G_S$ is a tree then
every edge appears twice.

\medskip
\noindent%
{\bf Corners of boundary cycles.} We have defined {\it corners} in
combinatorial pseudo-triangulations as being the angles marked
{\it small}. We extend the concept to the {\it vertices} on
boundary cycles of induced subgraphs $G_S$ by looking at the
labels of angles in $G$ incident to $v$ on the outer face of
$G_S$. We call $v$ a {\it corner of type $1$} if it contains a
{\it big} label on the outer face, or a {\it corner of type $2$}
when $v$ is non-pointed in $G$ but contains two consecutive {\it
small} labels on the outer face.

\smallskip
The following simple counting lemma will be useful later.

\begin{lemma}
\label{lemma:cornercount}%
Let $G_S$ be a subgraph of a cpt induced by the subset $S\subset
V$. Assume that $G_S$ is connected and that it contains all the
edges lying in the interior of its boundary cycle.

Let $m$ be the number of edges, $k$ the number of pointed
vertices, $l$ the number of non-pointed vertices in $G_S$ and $b$
the length of the boundary cycle in $G_S$. Then the number $c_1$
of corners of the type $1$ (big angles in the outer boundary) of
$G_S$ equals
\[
c_1 = m + 3 - 2k - 3l + b.
\]
\end{lemma}

In this statement a vertex in $G_S$ is called {\it pointed} if and
only if it was pointed in $G$.

\medskip
\begin{pf}
Let $f$ is the number of interior pseudo-triangles. The number of
interior angles in $G_S$ is $3f + k - c_1$, because there are $3f$
small interior angles and $k-c_1$ interior big angles. But the
number of interior angles also equals $2m-b$ (since the total
number of angles in any plane graph equals $2m$). Hence,
$2m-b=3f+k-c_1$, or
\[
3(m-f)=m+k +b-c_1.
\]
Finally,  Euler's formula applied to $G_S$ (as it contains all its
interior edges) is  $m-f=(k+l)-1$, which implies $3k+3l-3=m +k
+b-c_1$ and thus the desired statement.
\end{pf}

\medskip
\noindent {\bf The partially directed  auxiliary graph $D$ of a
combinatorial pseudo-triangulation $G$.} A \emph{partially
directed graph} $D=(V,E,\vec{E})$ is a graph $(V,E)$ together with
an assignment of directions to some of its edges. Thus edges are
allowed to get two directions, one direction only, or remain
undirected. Formally, $\vec{E}$ is a subset of $E\cup (-E)$ where
$E$ is the set of directed edges of $G$.

A plane embedding of a partially directed graph $(V,E,\vec{E})$ is
\emph{3-connected to the boundary} if from every interior vertex
$p$ there are at least three vertex-disjoint directed paths in
$\vec{E}$ ending in three different boundary vertices.
Equivalently, if for any interior vertex $p$ and for any pair of
forbidden vertices $q$ and $r$ there is a directed path from $p$
to the boundary not passing through $q$ or $r$.

\smallskip

\begin{lemma}
\label{lemma:Gconstruction} For every combinatorial
pseudo-triangulation, there exists a partially directed graph $D$
satisfying the following conditions:
\begin{enumerate}
\item $D$ is planar and contains the underlying graph of $G$. %
\item Every interior pointed vertex $v\in V$  of $G$ has three
out-neighbors: its two neighbors in $G$ along extreme edges and a
neighbor along the interior of the pseudo-triangle containing the
big angle at $v$. %
\item For every pointed vertex $v$ of $G$ its out-neighbors in $D$
are exactly its neighbors in $G$.
\end{enumerate}
\end{lemma}

\begin{figure}[ht]
\begin{center}
\ \psfig{figure=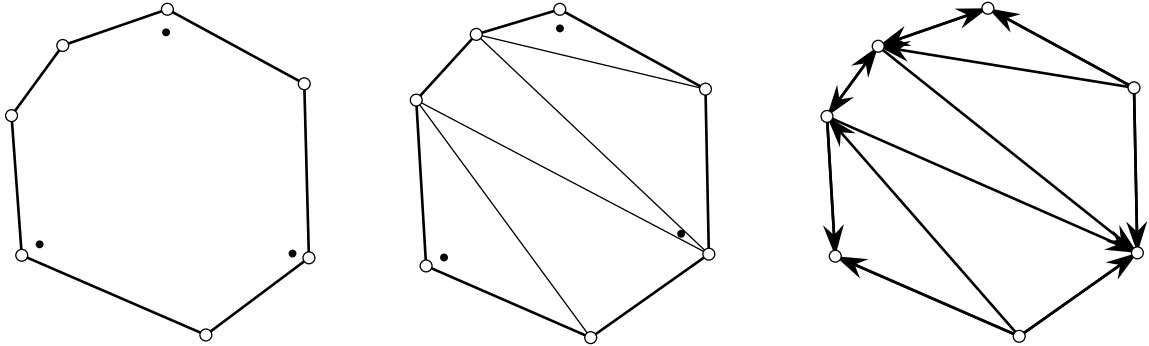,height=2.0in}
\end{center}
\caption{Left: A combinatorial pseudo-triangular face, with a
small black dot indicating a small angle (big angles are not
marked). Middle: a compatible triangulation of the face. Right:
the edges of the auxiliary directed graph.}
\label{figure:Gconstruction}
\end{figure}

\begin{pf}
We extend the underlying graph of $G$ to a (topological)
triangulation by triangulating the pseudo-triangles of $G$ with
more than three vertices in such a way that every big angle of $G$
is dissected by at least one new edge. This can be achieved by
recursively dissecting each face with an edge joining a
non-pointed vertex on the face to the opposite corner. Then the
edges are oriented as required by the statement.
See Figure \ref{figure:Gconstruction} for an illustration of how a
face is triangulated and how the edges incident to big angles are
oriented.
\end{pf}

\medskip

The main result of this section can now be stated.

\begin{theorem}
\label{thm:stretchable}%
For a combinatorial pseudo-triangulation $G$ with non-degenerate
(simple polygonal) faces the following are equivalent:
\begin{enumerate}
\item $G$ can be stretched into a compatible pseudo-triangulation.%
\item Every subgraph of $G$ with at least three vertices has at
least three corners. %
\item Every partially directed graph satisfying the requirements
of Lemma \ref{lemma:Gconstruction} is $3$-connected to the
boundary.
\end{enumerate}
\end{theorem}

\subsection{Proof of Theorem \ref{thm:stretchable}}

The implication from part 1 to part 2 is trivial. If  $G$ is
embedded as a pseudo-triangulation, there is no loss of generality
in assuming that the embedding is in general position, so that
every subgraph with at least three vertices has at least three
convex hull vertices. And all convex hull vertices of a subgraph
of $G$ will be corners according to our definition.

\medskip
\noindent {\bf Proof of $2\Rightarrow3$ in Theorem
\ref{thm:stretchable}:}\ \ To prove $3$-connectedness, we will
show that from every interior vertex $a$ there is a directed path
in $D$ going to the boundary and not passing through two arbitrary
(but fixed) vertices $b$ and $c$.

Let us consider $A$, the {\em directed connected component} of
vertex $a$, defined as the set of all vertices and directed edges
of  $D$ that can be reached from $v$ without passing through $b$
or $c$. We define this component as not containing the forbidden
points $b$ and $c$, but it may contain edges {\em arriving at
them}. Our goal is to prove that $A$ contains a vertex on the
boundary of $D$. We argue by contradiction. Suppose that all the
vertices of $A$ are interior to $D$.

For each interior pointed vertex $v$, let  $T_v$ be the unique
pseudo-triangle of $G$ containing the big angle at $v$. Thus $v$
is in an edge-chain of $T_v$ containing also the two extreme
adjacent edges of $v$. Let $G_S=(S,E_S)$ be the graph enclosing
all the pseudo-triangles $T_v$ associated to the pointed vertices
$v$ of $A$ and all the pseudo-triangles incident to the
non-pointed vertices. Clearly $G_S$ contains $A$: indeed, every
directed edge of $G$ is contained in a pseudo-triangle associated
to its source vertex.

We now use our hypothesis that $G_S$ has at least three corners.
We claim that at least one of them, $d$, belongs to $A$. This
gives the contradiction, because then there is an edge of
$D\setminus G$ jumping out of that corner $d$ (by the conditions
imposed to the partial orientation in $D$), which means that the
pseudo-triangle(s) corresponding to $d$ should have been contained
in $G_S$ and hence $d$ is not a corner of $G_S$ anymore.

To prove the claim, let $v_1,\dots,v_k$ be the corners of $G_S$
which are not in $A$. We want to prove that $k\le 2$. For this let
$T_1,\dots,T_k$ be pseudo-triangles in $G_S$, each having the
corresponding $v_i$ as a corner (there may be more than one  valid
choice  of the $T_i$'s; we just choose one). By definition,
some non-corner pointed vertex or some corner non-pointed vertex
of each $T_i$ is in the component $A$. Were there no forbidden
points, from a non-corner pointed vertex we could arrive to the
three corners of $T_i$ by three disjoint paths: two of them along
the concave chain containing the initial point and the third
starting with an edge of $D\setminus G$. From a non-pointed corner
vertex we could arrive to the other two corners by two disjoint
paths: moving out from the vertex to the two neighbors in the
incident pseudo-edges and then following along them. In
particular, since $v_i$ is not in $A$, one of the forbidden points
must obstruct to one of these paths, which implies that  either
$v_i$ equals one of the forbidden points $b$ and $c$ or that $T_i$
is the pseudo-triangle of one of the two forbidden points. And,
clearly, each of the two forbidden points contributes to only one
of the indices $i$ (either as a corner of $G_S$ or via its
associated pseudo-triangle if it is not a corner, but not both).
This shows that $k\le 2$ and completes the lemma.
\qed
\bigskip

\noindent{\bf Tutte's equilibrium condition.} To prove
$3\Rightarrow 1$ we use a directed version of Tutte's Theorem on
barycentric embeddings of graphs.

An embedding $D(P)$ of a partially directed graph
$D=(V,E,\vec{E})$ on a set of points $P=\{p_1, \cdots, p_n\}$,
together with an assignment $w:\vec{E}\rightarrow\R$ of weights to
the directed edges is said to be  in {\em equilibrium} at a vertex
$i\in V$ if
\[
\sum_{(i,j)\in\vec{E}} w_{ij} (p_i-p_j)=0.
\]

\begin{theorem} {\bf (Directed Tutte Theorem)}
\label{thm:directedTutte} Let $D=(\{1,\ldots,n\},E,\vec{E})$ be a
partially directed plane graph, $3$-connected to the boundary and
whose boundary cycle has no repeated vertices. Let
$(k+1,\ldots,n)$ be the ordered sequence of vertices in this
boundary cycle and let $p_{k+1},\ldots,p_n$ be the ordered
vertices of a convex $(n-k)$-gon. Let be the ordered vertices of a
convex $(n-k)$-gon. Let $w:\vec{E'}\rightarrow\R$ be an assignment
of positive weights to the internal directed edges. Then:
\begin{itemize}
\item[(i)] There are unique positions $p_1,\ldots,p_k\in\R^2$ for
the interior vertices such that all of them are in equilibrium in
the embedding $D(P), \ P=\{p_1,\ldots,p_n\}$.

\item[(ii)] In this embedding, all cells of $D$ are then realized
as non-overlapping convex polygons.
\end{itemize}
\end{theorem}

\begin{pf}
The proof of Tutte's Theorem given in \cite{richter-gebert}
(Theorem~12.2.2, pages~123--132) works with only minor
modifications. First, in the definition of {\em good
representation} (Definition~12.2.6, page~126), each point $p_i$ is
required to be in the relative interior of its {\em
out-neighbors}, since this is what the directed equilibrium
condition gives. Second, Claim~1 on page~126 proves that in a good
representation it is not possible for a vertex $p$ that $p$ and
all its neighbors lie in a certain line $\ell$, using
3-connectedness. The proof can be adapted to use {\em
$3$-connectedness to the boundary} as follows: consider three
vertex disjoint paths from $p$ to the boundary. Call $q$ any of
the three end-points, assumed not to lie in the line $\ell$.
Complete the other two paths to end at $q$ using boundary edges in
opposite directions. This produces three vertex-disjoint paths
from $p$ to a vertex $q$ not lying on $\ell$. The rest needs no
change.
\end{pf}
\bigskip

\noindent {\bf Proof of $3\Rightarrow1$ in Theorem
\ref{thm:stretchable}:}\ \ Construct an auxiliary partially
directed graph $D$ in the conditions of Lemma
\ref{lemma:Gconstruction}, choose arbitrary positive weights for
its directed edges, and apply the Directed Tutte Theorem to it.
Since all weights are positive, the equilibrium condition on an
interior vertex $p$, together with  the convexity of faces that
comes from Tutte's theorem implies that every interior vertex is
in the relative interior of the convex hull of its out-neighbors.
The conditions on $D$ then imply that the straight-line embeddding
of $G$ so obtained has big and small angles distributed as
desired. \qed

\medskip
\noindent {\bf Time Analysis.} Suppose that we are given a cpt
that can be stretched. Tutte's theorem actually gives an algorithm
to find a stretching: construct the auxiliary graph $D$ of Lemma
\ref{lemma:Gconstruction}, choose coordinates for the boundary
cycle in convex position and arbitrary positive weights for the
directed edges, and then compute the equilibrium positions.

Everything  can be done in linear time, except for the computation
of the equilibrium position for the interior vertices. In this
computation one writes a linear equation for each interior vertex,
which says that the position of the vertex is the average of its
(out-)neighbors. The position of the boundary vertices is fixed.
It has been observed~\cite[Section 3.4]{cgt-cdgtt-96} that the
planar {\it structure} of this system of equations allows it to be
solved in $O(n^{3/2})$ time, using the $\sqrt n$-separator theorem
for planar graphs in connection with the method of Generalized
Nested Dissection (see~\cite{lrt-gnd-79,lt-apst-80} or
\cite[Section 2.1.3.4]{r-tscts-88}), or even in time $O(M(\sqrt
n))$, where $M(n)=O(n^{2.375})$ is the time to multiply two
$n\times n$ matrices.

\subsection{Laman and circuit combinatorial pseudo-triangulations can be stretched}

Not all combinatorial pseudo-triangulations can be stretched: see
for instance the first example in Figure \figref{cptN}. Its non
stretchability can be proved either by showing that the graph is
not Laman (while the graph of every pointed pseudo-triangulation
must be so) or by applying the characterization given in Theorem
\ref{thm:stretchable} and finding a subgraph with less than three
corners.

Our next goal is to prove that if the underlying graph of a
combinatorial pseudo-triangulation is Laman (or is a rigidity
circuit)
then it can be stretched. The proof uses the Laman counting
condition to show that every subgraph has at least three corners.
We recall that both Laman graphs and circuits have the property
that a subgraph induced on a subset of $k\ge 2$ vertices has at
most $2k-3$ edges. By taking the complementary set of vertices
(and edges) this is equivalent to:
\begin{itemize}
\item In a Laman graph with $n$ vertices, every subset of $k\le
n-2$ vertices is incident to at least $2k$ edges. %
\item In a rigidity circuit graph with $n$ vertices, every subset
of $k\le n-2$ vertices is incident to at least $2k+1$ edges.
\end{itemize}

\begin{theorem}
\label{thm:Lamanstretch} Every subgraph $G_S$  of a Laman
combinatorial pseudo-triangulation $G$ has at least $3$ corners.
Therefore $G$ can be stretched.
\end{theorem}

\begin{pf}
We show first that there is no loss of generality in assuming that
$G_S$ is {\em simply connected} (i.e. it is connected and contains
all the edges of $G$ interior to its contour) and that no edge
appears twice in the boundary cycle.

If $G_S$ has an edge which  appears twice on the boundary cycle,
its removal does not change the number of corners; indeed, each
end-point of such an edge that is a corner after the removal must
be a corner before as well. If $G_S$ is not connected, either some
connected component has at least three vertices or all the
vertices of $G_S$ are corners. If $G_S$ is connected but not
simply connected, then adding to $G_S$ the pseudo-triangles, edges
and vertices of $G$ that fill in the holes does not change the
number of corners.

We now observe that since $G$ is pointed, the equation in Lemma
\ref{lemma:cornercount} becomes
\[e=2k-3-(b-c),
\]
where $e$, $k$, $b$ and $c$ are the numbers of edges, vertices,
boundary edges and corners of $S$.

Let
$b_0$ be the number of boundary vertices of $G_S$ (which may be
smaller than $b$, if a boundary vertex appears twice in the
boundary cycle). We now consider the set of edges incident to
vertices in the interior of $G_S$. Since there are $k-b_0$
interior vertices, the (rephrased) Laman property tells us that
there are at least $2(k-b_0)$ such edges. On the other hand, these
edges are all interior to $G_S$, and the total number of interior
edges in $G_S$ is $e-k$. Hence:
\[
e-k = 2k + c - 3 -2 \ge 2k-2b_0
\]
which implies the desired relation $c \ge 3 + 2 b-2 b_0 \ge 3$.
\end{pf}

\begin{theorem}
\label{thm:circuitstretch} Every subgraph $G_S$  of a rigidity
circuit cpt $G$ has at least $3$ corners. Hence, $G$ can be
stretched.
\end{theorem}

\begin{pf}
As in Theorem \ref{thm:Lamanstretch} we may assume without loss of
generality that $G_S$ is connected, contains all the edges of $G$
enclosed by its boundary cycle and its boundary cycle has no
repeated edges.

Let $k$, $l$, $e$, and $b$ be the numbers of pointed vertices,
non-pointed vertices, edges and boundary edges of $S$,
respectively. By  Lemma \ref{lemma:cornercount} we have
\[
e= 2k + 3l + c_1 - b- 3.
\]

If $G_S$ has no interior edge then the statement is trivial:
either $G_S$ contains no pseudo-triangle and then all its vertices
are corners, or it contains only one pseudo-triangle and the three
corners of it are corners of $G_S$. If $G_S$ has at least one
interior edge then:
\[
e \ge 2k+2l-b +1.
\]
Indeed, let $A\subset S$ be the set of vertices interior to $G_S$,
so that its cardinality equals $k+l-b_0$, where $b_0\le b$ is the
number of boundary vertices in $G_S$. If $A$ is empty then
$b_0=k+l\le b$ and our inequality becomes $e\ge b+1 - 2(b-b_0)$,
which holds by the existence of at least one interior edge. If $A$
is not empty we apply the Laman condition to $G_A$, which says
that the number of interior edges of $G_S$ is at least
$2(k+l-b_0)+1$, hence the total number of edges in $G_S$  is at
least $2(k+l-b_0)+1+b = 2k+2l-b+1 + 2(b-b_0) \ge 2k+2l-b +1$.

The two formulas above imply that
\[
l + c_1 \ge 4
\]
and, since $l\le 1$ (because there is only one non-pointed vertex
in $G$), $c_1 \ge 3$.
\end{pf}

\medskip

This completes the proof of Theorem \ref{thm:mainC}.


\section{Conclusions and Open Problems}
\label{section:open}

We have shown that any combinatorial pseudo-triangulation of a
plane Laman graph or of a plane rigidity circuit is stretchable.
In this latter case, we may even prescribe the non-pointed vertex.
In addition, Laman-plus-one and Laman-plus-two graphs are also
stretchable, although we may not in general be able to prescribe
the outer face or the non-pointed vertices.

The Main Result stated in the Introduction has thus been extended
along several lines, leading to interesting combinatorial objects
to study and several open questions, some solved in this paper,
some left for the future. We end with a listing of the main
directions for further investigations.

\medskip
\noindent%
{\bf Embeddability of planar generically rigid graphs as
pseudo-triangulations.} The goal here is the clarification of the
connection between minimum (pointed) pseudo - triangulations of a
planar point set and triangulations (maximal planar graphs
embedded on the same point set). Triangles are pseudo-triangles,
and every triangulation is a pseudo-triangulation, but some or all
of the vertices of the embedding may not be pointed. All planar
graphs containing a Laman graph are rigid (although not minimally
so). Stratifying by the number of additional edges (besides a
minimally rigid substructure) added to a Laman graph, we want to
investigate realizability as triangulations with some prescribed
number of non-pointed vertices. In this paper, we solved the case
of one additional edge (via the special case of rigidity
circuits). We leave open the question of completing the
characterization for the whole hierarchy. Such an investigation
will shed light into new intrinsic properties of planar
triangulations, some of the best studied and still elusive objects
in Combinatorial Geometry. We make the following conjecture.

\begin{conj} \label{conj:equivalences}
Given a plane graph $G$, the following conditions are equivalent:
\begin{itemize}
\item[(i)] $G$ is generically rigid%
\item[(ii)] $G$ can be straightened as a pseudo-triangulation.
\end{itemize}
\end{conj}

\medskip
\noindent%
{\bf Combinatorial pseudo-triangulations and embeddings on
oriented matroids (pseudolines).} We have seen that not all planar
graphs admitting combinatorial pointed pseudo-triangular labelings
are Laman graphs. But those which are Laman also have
straight-line realizations. A further direction of research
emerging from our work is to study the connection between the
combinatorial pseudo-triangulations and realizations in the
oriented matroid sense (on pseudo configurations of points).

\medskip
\noindent%
{\bf Grid size of pseudo-triangular embeddings.} Every planar
graph can be embedded on a grid of size $O(n)\times O(n)$, see for
example \cite{fpp}, \cite{schnyder}, \cite{felsner}. Here is a
natural remaining problem.

\smallskip
\noindent
\begin{open}
Can a planar Laman graph be embedded as a pseudo triangulation on
a $O(n^k) \times O(n^k)$ size grid? What is the smallest such $k$?
\end{open}

\medskip
\noindent%
{\bf Reciprocal duals of pseudo-triangulations.} Planar graphs
have combinatorial duals, obtained by replacing faces with
vertices and vice-versa. Moreover, when an embedded planar graph
supports a self-stress, it has a {\it geometric dual}, the
so-called {\it reciprocal diagram} of Maxwell \cite{
crapowhiteley}. A natural question (which will be answered in a
subsequent paper) concerns the connection between stressed
pseudo-triangulations (necessarily not minimal) and the planarity
of their reciprocal duals, see \cite{recip-duals}.

\medskip
\noindent%
{\bf Algorithmic questions.} We conjecture that our embedding
algorithms can be improved from $O(n^{3\over 2})$ to $O(n \log n)$
time. In general, the time complexity would be improved by a
positive answer to the following open questions.

\begin{open}
\label{open:LamanPlanar} Is it possible to decide the Laman
condition in linear time for a planar graph?
\end{open}

\begin{open}
\label{open:edgeBack} Is it possible to decide, faster than by
testing the Laman condition, which edge to put back in a Henneberg
II step for a planar graph? For a combinatorial
pseudo-triangulation?
\end{open}

\medskip
\noindent {\bf Acknowledgements} This research was initiated at
the {\it Workshop on Rigidity Theory and Scene Analysis} organized
by Ileana Streinu at the Bellairs Research Institute of McGill
University in Barbados,  Jan. 11-18, 2002 and partially supported
by NSF grant CCR-0203224. A preliminary version of this paper
appeared in \cite{socg03}.




\end{document}